\documentclass{article}
\usepackage{graphicx} 
\usepackage{amsmath}
\usepackage{comment}
\usepackage[shortlabels]{enumitem}
\usepackage[export]{adjustbox}
\usepackage[style=alphabetic, maxbibnames=99]{biblatex}
\AtEveryBibitem{
\clearfield{issn}
  \clearfield{isbn}
  \clearlist{location}
  \clearfield{month}
  \clearfield{day}
  \clearlist{language}
  \clearlist{translator}
  \renewbibmacro{in:}{}
}

\usepackage{subcaption}
\usepackage[shortlabels]{enumitem}

\newcounter{proofstep}
\newcommand{\step}{
  \refstepcounter{proofstep}
  \medskip
  \textit{Step \arabic{proofstep}.}
}
\AtBeginEnvironment{proof}{\setcounter{proofstep}{0}}

\newcounter{paragraphcount}[section]
\newcommand{\titledparagraph}[1]
{\bigskip\noindent%
\refstepcounter{paragraphcount}%
\textbf{\thesection.\theparagraphcount.\ #1}}

\newcommand{\titledparagraphnonr}[1]
{\bigskip\noindent%
\textbf{#1}}

\usepackage{fullpage}
\usepackage{mathtools}

\usepackage{amsmath}
\usepackage{amsthm}
\usepackage{thmtools}
\usepackage{hyperref}
\usepackage{tikz}
\usetikzlibrary{arrows.meta}
\usepackage{amssymb, amsfonts}

\usepackage{tikz-cd}
\usepackage{enumitem}

\usepackage{mathdots}

\usepackage{url}

\usepackage{wrapfig}

\usepackage[T1]{fontenc}

\usepackage{floatflt}

\theoremstyle{plain}
\newtheorem{theorem}{Theorem}[section]
\newtheorem{lemma}[theorem]{Lemma}
\newtheorem{proposition}[theorem]{Proposition}
\newtheorem{corollary}[theorem]{Corollary}

\newtheorem{question}[theorem]{Question}

\theoremstyle{definition}

\newtheorem{ex}[theorem]{Example}
\theoremstyle{remark}
\newtheorem{remark}[theorem]{Remark}

\usepackage{cleveref}

\newcommand{\C}{\mathbb{C}}
\newcommand{\R}{\mathbb{R}}
\newcommand{\Q}{\mathbb{Q}}
\newcommand{\Z}{\mathbb{Z}}
\newcommand{\p}{\mathbb{P}}

\newcommand{\Mod}{\mathrm{Mod}}
\newcommand{\MW}{\mathrm{MW}}
\newcommand{\MWH}{\mathrm{MW}^\mathrm{hol}}
\newcommand{\NS}{\mathrm{NS}}

\newcommand{\triv}{\mathrm{Triv}}
\newcommand{\diff}{\mathrm{Diff}}
\newcommand{\PD}{\mathrm{PD}}
\newcommand{\SL}{\mathrm{SL}}

\newcommand{\id}{\mathrm{id}}
\newcommand{\incl}{\mathrm{incl}}

\title{A Smooth Shioda--Tate Isomorphism for Elliptically Fibered 4-Manifolds}
\author{Maria Morariu}

\begin{document}

\maketitle

\begin{abstract}
    This paper studies a subgroup of the mapping class group of an elliptically fibered $4$-manifold ${\pi\colon M\to B}$, inspired by number theory. The smooth Mordell--Weil group $\MW(\pi)$ is the group of isotopy classes of smooth sections of $\pi$. We construct a smooth analog of the Shioda--Tate isomorphism, identifying $\MW(\pi)$ with the quotient of $H^2(M;\Z)$ by the trivial lattice of $\pi$.
    This shows that $\MW(\pi)$ can be computed explicitly from $H^2(M;\Z)$. As an application, we show that the natural inclusion of the holomorphic Mordell--Weil group into the smooth one is an isomorphism if and only if $M$ is a rational elliptic surface. 
    As another application, we show that $\MW(\pi)$ satisfies a strong version of Nielsen realization: any torsion element of $\MW(\pi)$ is represented by a holomorphic translation of $M$ of the same order.

\end{abstract}

\section{Introduction}
We start by giving our main theorem and then explain the objects involved in it. 
For an elliptic fibration $\pi\colon M\to B$ with fixed zero-section $\sigma_0$, the \textit{smooth Mordell--Weil group} is 
\[
    \MW(\pi)=\pi_0(\{\text{smooth sections of }\pi\}).
\]
The following is a smooth analog of the classical Shioda--Tate isomorphism between the holomorphic Mordell--Weil group $\MWH(\pi)$ of an elliptic surface and the quotient of its Néron--Severi lattice by the trivial lattice (\cite[Theorem 6.5]{Schutt_Schioda_MordellWeil}).

\begin{theorem}\label{thm: full MW} Let $\pi\colon M\to B$ be an elliptic fibration over a closed base curve $B$ with fixed zero-section $\sigma_0\colon B\to M$. Then there exists a natural group isomorphism
    \[
        \MW(\pi) \xrightarrow[]{\cong} H^2(M;\Z)/\triv(\pi)
    \]
    \[
        [\sigma]\mapsto [\PD_M(\sigma(B))-\PD_M(\sigma_0(B))],
    \]
where $\sigma(B)$ is viewed as a homology class in $H_2(M;\Z)$ and $\PD_M(\sigma(B))$ denotes its Poincaré dual.
\end{theorem}

We apply \Cref{thm: full MW} in concrete examples below. One useful consequence is that for any elliptic surface over a closed base, $H^2(M)\coloneqq H^2(M;\Z)$ is generated by horizontal elements (elements corresponding to sections) and vertical elements (ones of the trivial lattice). As a key application of \Cref{thm: full MW} and the Shioda--Tate isomorphism, we show the following corollary.

\begin{corollary}[\textbf{Characterizing rational elliptic surfaces}]\label{cor: rat ell surf characterization}
    For an elliptic surface $\pi\colon M\to B$ over a closed base curve $B$ with holomorphic zero-section $\sigma_0$, the groups $\MW(\pi)$ and $\MWH(\pi)$ are isomorphic if and only if $M$ is rational.
\end{corollary}

\titledparagraph{Main Objects.}
We introduce the main objects of this paper.

An \textit{elliptic fibration} is a triple $(M,B,\pi)$, where $M$ is a smooth complex surface, $B$ is a smooth complex curve and $\pi\colon M\to B$ is a proper, surjective, holomorphic map whose general fiber is a smooth genus $1$ curve. 
In this paper, we further assume that an elliptic fibration comes equipped with a smooth zero-section $\sigma_0\colon B\to M$, that $\pi$ has at least one critical point, and that $\pi$ is relatively minimal (i.e.\ no fiber contains a complex curve of self-intersection $-1$). The choice of zero-section equips each smooth fiber with a distinguished point, which gives it an abelian group structure. This group structure extends to the smooth locus of the singular fibers, as explained in \Cref{section: group structure fibers}.

Let $\Gamma^\mathrm{sm}(\pi)$ denote the group of smooth sections of $\pi\colon M\to B$, where sections are added pointwise and the identity is the distinguished section $\sigma_0$.
The \textit{(smooth) Mordell--Weil group} of an elliptic fibration $\pi\colon M\to B$ is
\[
    \MW(\pi)\coloneqq \pi_0(\Gamma^\text{sm}(\pi)).
\]
Replacing smooth sections by holomorphic sections of $\pi$ and stipulating that the zero-section is holomorphic yields the classical (holomorphic) Mordell--Weil group $\MWH(\pi)$. If the base $B$ is closed, both $\MW(\pi)$ and $\MWH(\pi)$ are finitely generated abelian groups (see \cite[Theorem 3.4]{FarbLoojenga_MordellWeil} and \cite[Theorem 6.6]{Schutt_Schioda_MordellWeil}).

The \textit{trivial lattice} associated to an elliptic fibration $\pi\colon M\to B$ is the subgroup of $H^2(M)$ generated by the Poincaré duals of the zero-section and the irreducible components of all fibers of $\pi$.

Using \Cref{thm: full MW}, we compute the Mordell--Weil group of two elliptic K3 surfaces. We fill in the details of these examples in \Cref{subsection: applications}.

\begin{figure}
    \centering
    \begin{subfigure}[b]{0.45\textwidth}
        \centering
        \includegraphics[width=\textwidth]{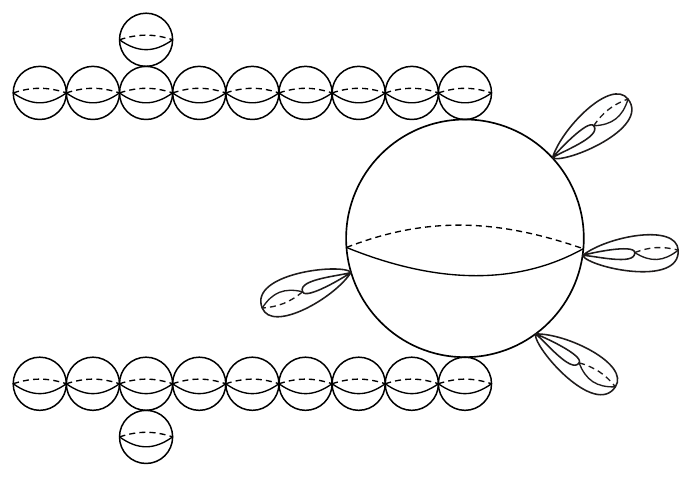}
        \caption{The Inose K3 surface has $\MW(\pi)\cong \Z^4$.}
        \label{fig: inose intro}
    \end{subfigure}
    \hspace{1.5cm}
    \begin{subfigure}[b]{0.3\textwidth}
        \centering
        \includegraphics[width=\textwidth]{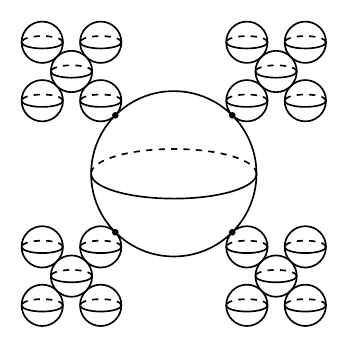}
        \caption{The Kummer K3 surface has $\MW(\pi)\cong \Z^4\times \Z/2\Z\times \Z/2\Z$.}
        \label{fig: kummer intro}
    \end{subfigure}
    \caption{Two elliptic K3 surfaces with different Mordell--Weil groups. In both figures, the large sphere represents a section of $\pi$. The smaller spheres and nodes are components of the singular fibers of the two fibrations. All fibers which are not depicted are smooth tori.}
\end{figure}

\begin{ex}[\textbf{Mordell--Weil Group of Inose K3 Surface}]\label{ex: intro Inose}
    Consider the Inose K3 surface with two singular fibers of type $II^*$ and four singular fibers of type $I_1$ as in \Cref{fig: inose intro}. The Poincaré duals of irreducible components of the two fibers of type $II^*$ which do not intersect the zero-section, the fiber class and the zero-section generate the trivial lattice. By computing intersection numbers as in Example \ref{ex: inose surface} in a later section, we show that the trivial lattice is unimodular and thus primitively embedded in $H^2(M)$. Using vanishing cycles and Lefschetz thimbles, we can give the remaining $4$ generators of $H^2(M)$. \Cref{thm: full MW} implies
    \[
        \MW(\pi)\cong \Z^4.
    \]
\end{ex}

\begin{ex}[\textbf{Mordell--Weil Group of Kummer Surface}]\label{ex: intro Kummer}
    Let $\pi\colon M\to \p^1$ be an elliptically fibered Kummer surface as in \Cref{fig: kummer intro}. It has $4$ singular fibers of type $I_0^*$. The trivial lattice is generated by the Poincaré duals of the irreducible components of the singular fiber, the class of the generic fiber and the zero-section. In \Cref{subsection: applications}, we give explicit generators for $H^2(M)$ and conclude using \Cref{thm: full MW} that 
    \[
        \MW(\pi)\cong \Z^4\times \Z/2\Z \times \Z/2\Z.
    \]
\end{ex}
The elliptic surfaces in Example \ref{ex: intro Inose} and Example \ref{ex: intro Kummer} are both K3 surfaces and thus diffeomorphic, as shown by Kodaira. The two examples show that the Mordell--Weil group detects the difference between the two elliptic fibration structures.

\begin{figure}
    \centering
    \begin{subfigure}[b]{0.3\textwidth}
        \centering
        \includegraphics[width=\textwidth]{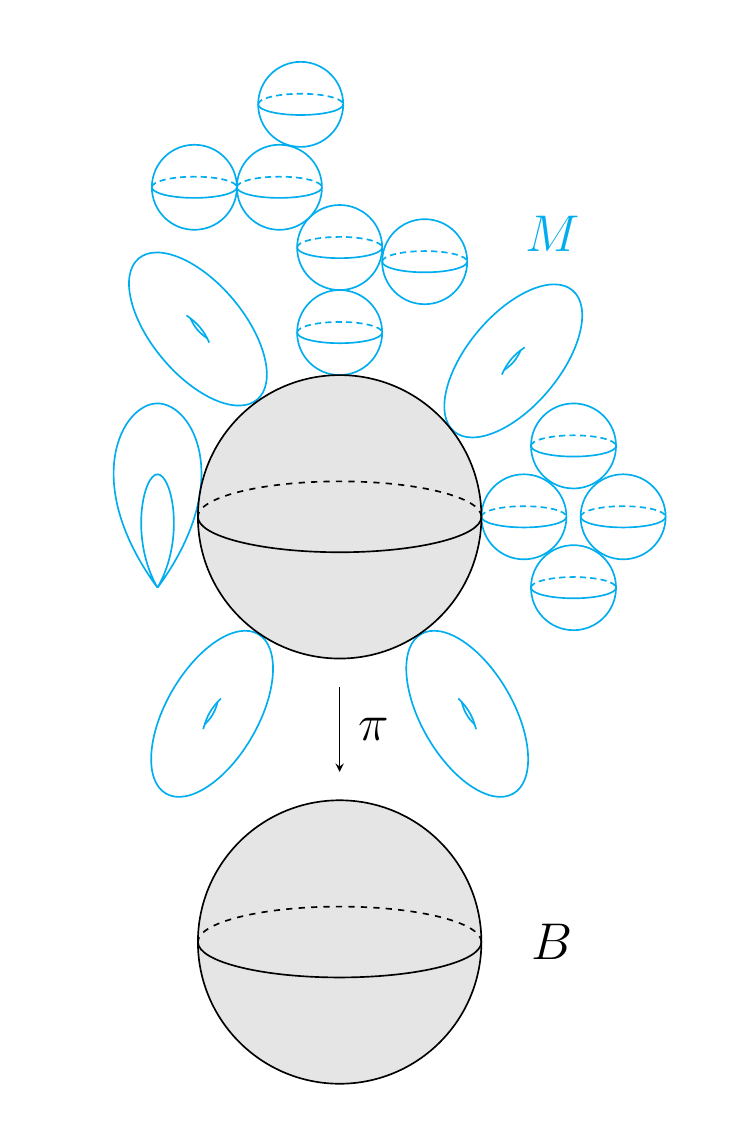}
        \caption{An elliptic fibration $\pi\colon M\to B$ with three singular fibers.}
    \end{subfigure}
    \hspace{1.5cm}
    \begin{subfigure}[b]{0.3\textwidth}
        \centering
        \includegraphics[width=\textwidth]{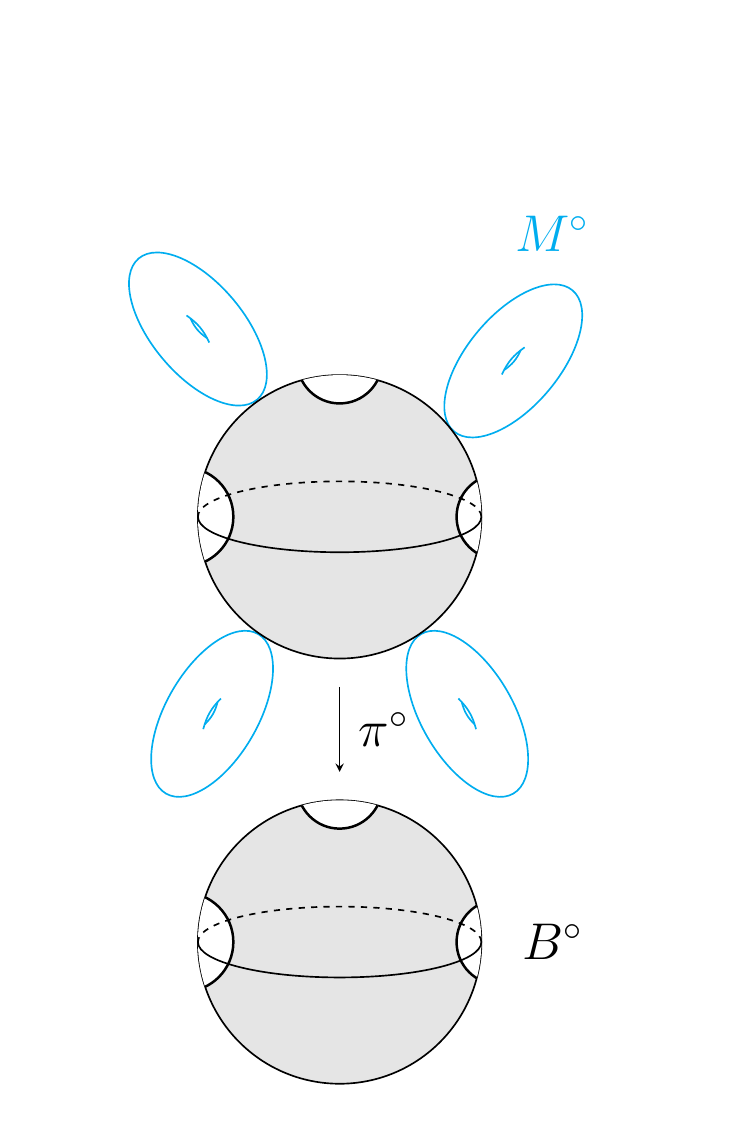}
        \caption{The reduced elliptic fibration $\pi^\circ\colon M^\circ\to B^\circ$ associated to $\pi$.}
    \end{subfigure}
    \caption{$M^\circ$ is obtained by removing an open neighborhood of each singular fiber of $M$.} 
    \label{fig: intro pi and pi circ}
\end{figure}

\titledparagraph{Proof Strategy.}
The main idea in proving \Cref{thm: full MW} is to give a cohomological characterization for $\MW(\pi)$ in terms of sheaf cohomology, similarly to \cite[Proposition 3.3]{FarbLoojenga_MordellWeil} and \cite[Proposition 3.7]{Morariu_NarrowMW}.
Such a characterization is made difficult by the possibility of a nontrivial component group of the relative smooth locus of reducible singular fibers, so we proceed by removing the singular fibers.

For an elliptic fibration $\pi\colon M\to B$, let $\Sigma\subset B$ be the set of critical values of $\pi$. For each $t\in \Sigma$, consider a small open disk $D_t$ around $t$ and assume that the disks $D_t$ are pairwise disjoint. Define
\[
    B^\circ\coloneqq B-\bigcup_{t\in \Sigma} D_t
\]
and the compact $4$-manifold with boundary 
\[
    M^\circ\coloneqq \pi^{-1}(B^\circ).
\]
Denote by $\pi^\circ\colon M^\circ\to B^\circ$ the restriction of $\pi$ (see \Cref{fig: intro pi and pi circ}).  
Any section $\sigma\colon B^\circ\to M^\circ$ naturally represents an element $[\sigma]\in H_2(M^\circ,\partial M^\circ)$.

\medskip\textit{Step 1. Construct an isomorphism between $\MW(\pi^\circ)$, the group of isotopy classes of sections of $\pi^\circ$, and the annihilator of the fiber class in $H^2(M^\circ)$.} We first express $\MW(\pi^\circ)$ in terms of sheaf cohomology (\Cref{prop: coh characterization of MW(pi^circ)}). Using a spectral sequence argument, similar to the one in \cite{FarbLoojenga_MordellWeil} and \cite{Morariu_NarrowMW}, we construct the desired isomorphism (\Cref{prop: MW bundle no singularities}).

\medskip\textit{Step 2. Relate the groups $\MW(\pi)$ and $\MW(\pi^\circ)$.} For $t\in B$, let $M_t\coloneqq \pi^{-1}(t)$. 
We construct a map $\Psi\colon \MW(\pi^\circ)\to \bigoplus_{t\in \Sigma} H_1(M_t)$, which can be described geometrically as in \Cref{section: globally: relation between MW(pi) and MW(pi^circ)}. The following theorem shows that $\Psi$ is the only obstruction to extending sections of $\pi^\circ$ to sections of $\pi$ up to isotopy.

\begin{theorem}[\textbf{Relation between $\MW(\pi)$ and $\MW(\pi^\circ)$}] \label{thm: relation MW pi and pi circ complete}
    For any elliptic fibration $\pi\colon M\to B$ over a closed base, there exists the following short exact sequence of groups:
    \[
    0\to \MW(\pi) \to \MW(\pi^\circ) \xrightarrow[]{\Psi} \bigoplus_{t\in \Sigma} H_1(M_t) \to 0.
    \]
\end{theorem}

\textit{Step 3. Conclusion} Using Step 1, we show that the kernel of $\Psi$ is isomorphic to $H^2(M)/\triv(\pi)$, which by Step 2 proves \Cref{thm: full MW}.

\titledparagraph{$\MW(\pi)$ as a Subgroup of the Mapping Class Group.} The Mordell--Weil group can be rephrased in terms of diffeomorphisms as
\[
    \MW(\pi)= \pi_0(\{f\in \diff(M)\mid f \text{ acts on each smooth fiber by translation}\}).
\]
Indeed, a diffeomorphism $f\in \diff(M)$ that acts by translation on each smooth fiber determines a smooth section $f\circ\sigma_0$ of $\pi$. Conversely, any smooth section $\sigma$ determines a unique diffeomorphism $f_\sigma\in \diff(M)$, that acts by translation by $\sigma$ on the smooth locus of $M$ (see \cite[Lemma 3.1]{Morariu_NarrowMW}).

A consequence of \Cref{thm: full MW} is that $\MW(\pi)$ acts faithfully on $H^2(M)$. 
This action factors as 
\[
    \MW(\pi) \to \pi_0(\diff(M)) \to \mathrm{Aut}(H^2(M),\smile)
\]
and injectivity of the composition implies that the map $\MW(\pi)\to \pi_0(\diff(M))$ is injective.

A strong version of Nielsen realization holds for Mordell--Weil elements: every torsion element of the Mordell--Weil group has a holomorphic representative that is torsion.
\begin{theorem}[\textbf{Nielsen realization}]\label{thm: nielsen realization}
    If $\pi\colon M\to B$ is an elliptic fibration over a closed base curve $B$ with a holomorphic zero-section $\sigma_0$, then
    the natural inclusion $\MWH(\pi)\xhookrightarrow{} \MW(\pi)$ induces an isomorphism of the torsion subgroups of these finitely generated abelian groups. 
\end{theorem}
\Cref{thm: nielsen realization} provides a practical method for computing $\MW(\pi)$. Let $g$ be the genus of $B$. Compute the number of additive and multiplicative singular fibers of $\pi$, denoted $a$ and $m$.
\cite[Corollary 1.4]{Morariu_NarrowMW} implies that the rank of $\MW(\pi)$ is $4g+2a+m-4$.
Then determine the torsion part using classical methods.

\titledparagraph{Related Work.} This paper extends previous work of the author \cite{Morariu_NarrowMW}. The previous paper computed the narrow Mordell--Weil group, a finite-index subgroup of $\MW(\pi)$ of an elliptic fibration. The smooth Mordell--Weil group was first introduced by Farb and Looijenga \cite{FarbLoojenga_MordellWeil}. They proved finite generation and computed the smooth Mordell--Weil group of a generic (i.e.\ assuming that all fibers are nodal) elliptic fibration $\pi\colon M\to \p^1$. A symplectic version of the Mordell--Weil group under the assumption that all singular fibers are nodal was computed by Hacking-Keating \cite{HackingKeating_Symplectomorphisms}. Related ideas were used by Cox and Zucker \cite{CoxZucker_IntersectionNumbers}, who used an injection of the tensor of holomorphic Mordell--Weil group with $\Q$ into $H^1(B,R^1\pi_*\Q)$ for computations.

\titledparagraph{Outline.} In \Cref{section: background}, we collect some background on elliptic surfaces that is used throughout the paper. In \Cref{section: compute MW pi circ}, we give a cohomological characterization for $\MW(\pi^\circ)$ and express this group as the annihilator of the fiber class. \Cref{section: local picture} contains an extensive analysis of an elliptic fibration over a disk with only one singular fiber at the origin, with explicit computations of $\MW(\pi)$ and $\MW(\pi^\circ)$. \Cref{section: globally: relation between MW(pi) and MW(pi^circ)} establishes a complete obstruction to extending isotopy classes of sections of $\pi^\circ$ to isotopy classes of sections of $\pi$ up to isotopy. In \Cref{section: proof main thm}, we prove \Cref{thm: full MW}. In \Cref{subsection: applications}, we apply \Cref{thm: full MW} in some concrete examples.

\titledparagraphnonr{Computational Resource Disclosure.} Gemini Pro 3.1 and Claude Opus 5.5 were used in the creation of pictures. GPT-5.6 Terra was involved in checking the proofs and literature search.
The only mathematical use of artificial intelligence in this paper was in Example \ref{ex: inose surface}, which was computed in conversation with GPT-6.1 Sol.

\titledparagraphnonr{Acknowledgements.} 
I am very grateful to my advisor Benson Farb for suggesting this project, for his constant support, and for extensive comments on earlier drafts. I am thankful to Eduard Looijenga for many useful discussions. I would like to thank Eduard Looijenga, Dan Margalit and Katherine Booth for helpful comments.

\section{Background on Elliptic Surfaces}
\label{section: background}

We record the background which is necessary for the paper.

\medskip\noindent\textbf{Notation and Conventions:} 
For an elliptic fibration $\pi\colon M\to B$, and $b\in B$, we denote $M_t\coloneqq \pi^{-1}(t)$.
All (co)homology is taken with $\Z$-coefficients, unless otherwise specified. For any closed $n$-manifold $X$, we denote by $\PD_X\colon H_k(X)\to H^{n-k}(X)$ the Poincaré duality isomorphism in $X$. For a compact $n$-manifold $X$ with boundary, we denote by $\PD_{[X,\partial X]}\colon H^k(X)\to H_{n-k}(X,\partial X)$ and $\PD_{[X,\partial X]}\colon H^k(X,\partial X)\to H_{n-k}(X)$  the isomorphisms given by Poincaré--Lefschetz duality. We assume that the zero-section $\sigma_0$ is holomorphic in a neighborhood of each critical value\footnote{This assumption is not restrictive, since for any critical value $t\in B$, the map $\pi$ is a holomorphic submersion at $\sigma_0(t)$, so $\sigma_0$ can be deformed to a section that is holomorphic in a neighborhood of $t$.}. The topology on $\diff(M)$ and $\Gamma^\mathrm{sm}(\pi)$ is the Whitney $C^\infty$ topology.

\begin{proposition}\label{prop: H_1 of elliptic fibration}
    If $\pi\colon M\to B$ is an elliptic fibration over a closed base $B$, then $\pi$ induces an isomorphism \[\pi_*\colon H_1(M)\xrightarrow[]{\cong} H_1(B).\]
\end{proposition}
\begin{proof}
    We assume that an elliptic fibration is non-trivial and has a section. In particular, an elliptic fibration has no multiple fibers. Hence, \cite[p. 157, Proposition 2.1]{FriedmanMorgan_Smooth4Mflds} implies that $\pi$ induces an isomorphism $\pi_*\colon \pi_1(M)\to \pi_1(B)$. Applying abelianization and the Hurewicz isomorphism on both sides yields the conclusion.
\end{proof}

\begin{proposition}[\protect{\cite[Proposition 5.4]{Honda_Viaklovski_Fibrations_6Sphere}}]\label{prop: noighborhood of singular fiber strong deformation retracts}
    If $\pi\colon M\to D$ is an elliptic fibration over a disk with only one critical value at the origin, then $M$ strong deformation retracts to $M_0$.
\end{proposition}

\subsection{Monodromy and Kodaira's Classification of Singular Fibers}
Let $\pi\colon M\to B$ be an elliptic surface with set of singular values $\Sigma\subset B$, and let $M_t$ be a singular fiber.
As a divisor, it can be written $M_t=\sum m_i C_i$, where $C_i$ are the distinct irreducible components of $M_t$ and $m_i\ge 1$ is the multiplicity of the component $C_i$. 
The definition of a section gives $\pi\circ\sigma_0=\id_B$ and by the chain rule, $d\pi_{\sigma_0(t)}\circ d\sigma_0(t) = \id_{T_b B}$. Hence $d\pi_{\sigma_0(t)}$ has full rank and $\sigma_0(t)$ lies in a multiplicity-one component of $M_t$.
In particular, there exists $i$ such that $m_i=1$.\footnote{The assumption that $\pi$ has a section is crucial; otherwise the fibration could have multiple fibers.}
To a reduced singular fiber $M_t$ we associate a \textit{dual graph}: a graph with vertices corresponding to $C_i$ and $C_i\cdot C_j$ edges connecting the vertices corresponding to $C_i$ and $C_j$, if $i\ne j$. 

Let us fix a regular value $b\in B$ and let
\[
    \rho\colon \pi_1(B-\Sigma, b)\to \Mod(M_b)
\]
be the \textit{monodromy representation} of $\pi$.
For a critical value $t\in B$, the monodromy around $t$ is $\rho([\gamma])\in \Mod(M_b)$, where $\gamma$ is a loop based at $b$, containing only the critical value $t$. Choosing a basis for $H_1(M_b)$ provides an identification $\Mod(M_b) \cong \SL(2,\Z)$.

The singular fibers of an elliptic surface were first classified by Kodaira \cite{Kodaira_CompAnalyticSurf2} and later by Néron \cite{Neron_ModelMin} and Tate \cite{Tate_AlgorithmSingFibers}.
Below, we give the classification of singular fibers from \cite[Theorem I.6.6, Table VI.2.1]{Miranda_BTES}, excluding the possibility of multiple fibers.

\begin{figure}[h]
    \centering
    \includegraphics[width=0.9\linewidth]{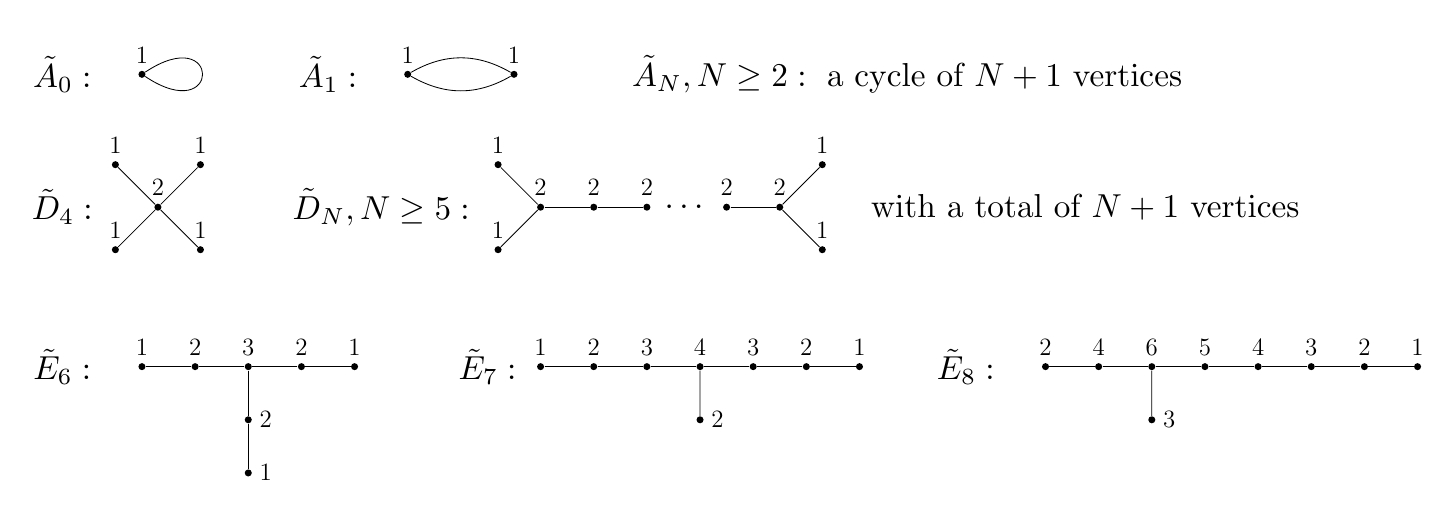}
    \caption{Extended Dynkin Diagrams of type ADE with component multiplicities.}
    \label{fig: extended Dynkin ADE}
\end{figure}

\begin{theorem}[\textbf{Kodaira Classification}]\label{thm: Kodaira classification sing fibers}
    The dual graph of a reducible fiber of an elliptic surface is an extended Dynkin diagram of type ADE and the multiplicities of the irreducible components are as shown in \Cref{fig: extended Dynkin ADE}.
    Furthermore, all fibers of an elliptic surface have one of the types in \Cref{table: kodaira classification}.
\end{theorem}
\begin{table}
    \centering
    \caption{Kodaira's classification of singular fibers of an elliptic surface.}\label{table: kodaira classification}
    \renewcommand{\arraystretch}{0.81}
    \begin{tabular}{c|c|c|c}
       \textbf{Kodaira Symbol}  & \textbf{Description} & \textbf{Dual Graph} & \begin{tabular}{c}
            \textbf{Monodromy up to}  \\
            \textbf{Conjugation} 
       \end{tabular} \\ \hline
       $I_0$  & smooth elliptic curve & &
       $\begin{pmatrix}
           1&0\\
           0&1
       \end{pmatrix}$\\ \hline
       $I_1$ & nodal rational curve & &
       $\begin{pmatrix}
           1&1\\
           0&1
       \end{pmatrix}$\\ \hline
       $I_N, N\ge 2$ &\begin{tabular}{c}
           $N$ smooth rational curves meeting \\
           transversally in a cycle 
       \end{tabular}
       & $\tilde{A}_{N-1}$&
       $\begin{pmatrix}
           1&n\\
           0&1
       \end{pmatrix}$\\ \hline
       $I_N^*, N\ge 0$ & $N+5$ smooth rational curves & $\tilde{D}_{N+4}$ &
       $\begin{pmatrix}
           -1&-n\\
           0&-1
       \end{pmatrix}$\\ \hline
       $II$ & a cuspidal rational curve & &
       $\begin{pmatrix}
           1&1\\
           -1&0
       \end{pmatrix}$\\ \hline
       $III$ & 
       \begin{tabular}{c}
            two smooth rational curves meeting  \\
            at one point with multiplicity $2$ 
       \end{tabular}
       & $\tilde{A}_{1}$ &
       $\begin{pmatrix}
           0&1\\
           -1&0
       \end{pmatrix}$\\ \hline
       $IV$ & 
       \begin{tabular}{c}
          three smooth rational curves meeting    \\
          at one point
       \end{tabular}
       & $\tilde{A}_{2}$ &
       $\begin{pmatrix}
           0&1\\
           -1&-1
       \end{pmatrix}$\\ \hline
       $IV^*$ & $7$ smooth rational curves& $\tilde{E}_6$ &
       $\begin{pmatrix}
           -1&-1\\
           1&0
       \end{pmatrix}$\\ \hline
       $III^*$ & $8$ smooth rational curves & $\tilde{E}_7$ &
       $\begin{pmatrix}
           0&-1\\
           1&0
       \end{pmatrix}$\\ \hline
       $II^*$ & $9$ smooth rational curves &  $\tilde{E}_8$ &
       $\begin{pmatrix}
           0&-1\\
           1&1
       \end{pmatrix}$
    \end{tabular}
    \end{table}

Consider a reducible singular fiber $M_t$ with irreducible components $C_0,\dots, C_r$, where $C_0$ intersects $\sigma_0$.
The dual graph of $M_t$ is an extended Dynkin diagram, so we can associate an \textit{extended Cartan matrix} $\tilde A$ with entries $\tilde A_{i,j}=-C_i\cdot C_j$. Let $\Bar{m}\in \Z^{r+1}$ be the vector of multiplicities of the components of $M_t$. Then the kernel of $\tilde A$ is spanned by $\Bar{m}$.
We can also associate a (non-extended) Dynkin diagram of type ADE by removing the vertex corresponding to $C_0$ from the dual graph of $M_t$. Thus, there exists a \textit{Cartan matrix} $A$, which satisfies
\[
    \tilde A=\begin{pmatrix}
        2 & v^T \\
        v & A
    \end{pmatrix}.
\]
The matrix $A$ has full rank and its determinant is equal to the number of components of multiplicity $1$ of $M_b$.

\subsection{The Smooth Locus of a Singular Fiber as a Group}\label{section: group structure fibers}
Consider an elliptic fibration $\pi\colon M\to B$.
The \textit{smooth locus} of $\pi$ is $M^\#\coloneqq \{x\in M: d\pi_x \text{ has full rank}\}$. For any $t\in B$, the smooth locus of the fiber $M_t$ is $M_t^\#\coloneqq M_t\cap M^\#$. If 
$M_t$ is smooth, $M_t^\#=M_t$. Otherwise, $M_t^\#$ is obtained by removing all singular points and higher-multiplicity components of $M_t$. 

An \textit{analytic fiber system of groups over a complex manifold $\Delta$} is a complex manifold $\mathcal{G}$ with a holomorphic submersion $\Psi\colon \mathcal{G}\to \Delta$ and a holomorphic section $e\colon \Delta\to \mathcal{G}$, such that
\begin{itemize}
    \item For all $t \in \Delta$, the preimage $\Psi^{-1}(t)$ is a Lie group with identity element $e(t)$.
    \item Denoting 
    \[
    \mathcal{D}\coloneqq \{(g_1,g_2)\in \mathcal{G}\times \mathcal{G}: \Psi(g_1)=\Psi(g_2)\}\subset \mathcal{G}\times \mathcal{G},
    \]
    the map $\mathcal{D}\to \mathcal{G}$ given by $(g_1,g_2)\mapsto g_1g_2^{-1}$ is holomorphic.
\end{itemize}

We show that locally, the smooth locus of an elliptic fibration is an analytic fiber system of groups. Let $t\in B$ be a critical value and let $U_t\subset B$ be a disk around $t$ such that $\sigma_0|_{U_t}$ is holomorphic. Denote $U_t'\coloneqq U_t-\{t\}$ and $M_{U_t}'\coloneqq \pi^{-1}(U_t')$. Then $\pi\colon M_{U_t}'\to U_t'$ is an analytic fiber system of groups with the holomorphic section $\sigma_0|_{U_t}$.
The following shows that the group structure extends to the smooth locus of the singular fibers.

\begin{theorem}[\protect{\cite[Theorem 9.1]{Kodaira_CompAnalyticSurf2}}]\label{thm: group on sing fibers}
    Let $\pi\colon M \to B$ be an elliptic fibration and let $t\in B$ be a critical value. For a disk $U_t\subset B$ containing $t$, such that $\sigma_0|_{U_t}$ is holomorphic, let $M_{U_t}\coloneqq \pi^{-1}(U_t)$. Then there exists on $M_{U_t}^\#$ a unique structure of analytic fiber system of abelian groups over $U_t$, which is an extension of the structure on $M_{U_t}'$ of analytic fiber system of abelian groups over $U_t'$.
\end{theorem}

Let $t\in B$ be a critical value of $\pi$ and denote by $C_0$ the irreducible component of $M_t$ intersecting $\sigma_0$. Let $C_0^\#\coloneqq C_0\cap M_t^\#$ be the connected component of the zero-element in $M_t^\#$ and define $G(M_t)\coloneqq \pi_0(M_t^\#)$. These fit into the following short exact sequence of topological groups:
\[
    0\to C_0^\#\to M_t^\# \to G(M_t) \to 0.
\]
Topologically, $C_0^\#$ is a sphere with one or two points removed (by Kodaira's classification). Hence, as a group, $C_0^\#$ is isomorphic to $\C$ or $\C^\times$. If $C_0^\#\cong \C$, the fiber $M_t$ is an \textit{additive fiber}, otherwise $M_t$ is a \textit{multiplicative fiber}. By Kodaira's classification, the only multiplicative fibers are those of type $I_n$ for $n\ge 0$.
The Lie group structure of $G(M_t)$, depending on the type of the singular fiber $M_t$ is given in \Cref{table: group structure on singular fibers}.
\begin{table}
    \centering
    \caption{The group structure on the smooth locus of $M_t$ is the semidirect product of $C_0^\#$ and $G(M_t)$.}\label{table: group structure on singular fibers}
    \renewcommand{\arraystretch}{1.3}
        \begin{tabular}{c|c}
        \textbf{Type of $M_t$} & \textbf{Group Structure of $G(M_t)$} \\ \hline
        $I_n,n\ge 1$ & $\Z/n\Z$ \\
        $I_{2n}^*$ & $(\Z/2\Z)^2$\\
        $I_{2n+1}^*$ & $\Z/4\Z$\\
        $II$ & $\{0\}$\\
        $III$ & $\Z/2\Z$\\
        $IV$ & $\Z/3\Z$\\
        $IV^*$ & $\Z/3\Z$\\
        $III^*$ & $\Z/2\Z$\\
        $II^*$ & $\{0\}$
    \end{tabular}
\end{table}

\begin{lemma} \label{lem: component group is monodromy coinvariants}
    Let $\pi\colon M\to B$ be an elliptic fibration. Fix a regular value $b\in B$ and a critical value $t\in B$. Let $f\in \Mod(M_b)$ be the monodromy around $M_t$ and let $f_*$ be the action it induces on $H_1(M_b)$. Then
    \[
        G(M_t)\cong (\Z^2/(I-f_*)\Z^2)_\text{tor}.
    \]
\end{lemma}
\begin{proof}
    This can be checked for each type of singular fibers using \Cref{thm: Kodaira classification sing fibers} and \Cref{thm: group on sing fibers}.
\end{proof}

\subsection{The Holomorphic Mordell--Weil Group}
The \textit{holomorphic Mordell--Weil group} of an elliptic surface $\pi\colon M\to B$ with holomorphic zero-section is the group of holomorphic sections of $\pi$. The correspondence between sections and diffeomorphisms that act by translation makes this equivalent to the group of biholomorphic automorphisms of $M$ that act by translation on each smooth fiber.

The \textit{Néron-Severi lattice} $\NS(M)$ of a complex surface $M$ is defined as the group of divisors of $M$, quotiented by algebraic equivalence. If the base curve $B$ is closed, the Lefschetz $(1,1)$-theorem implies
\[
    \NS(M)=H^2(M;\Z)\cap H^{1,1}(M).
\]

The following theorem gives an explicit identification of $\MWH(\pi)$ in terms of the Néron-Severi lattice and the trivial lattice. For a proof, see \cite[Theorem 6.5]{Schutt_Schioda_MordellWeil}.
\begin{theorem}[\textbf{Shioda--Tate Isomorphism}]
\label{thm: holomorphic Mordell--Weil} For any elliptic fibration $\pi\colon M\to B$ over closed base $B$ with holomorphic zero-section, there is an isomorphism
    \[
    \begin{array}{ccc}
         \MWH(\pi) &\xrightarrow[]{\cong}& \NS(M) / \triv(\pi),  \\
         \sigma & \mapsto &[\sigma].
    \end{array}
    \]
\end{theorem}

\section{Computing $\MW(\pi^\circ)$}
\label{section: compute MW pi circ}
Our computation of $\MW(\pi^\circ)$ relies heavily on the methods from \cite{FarbLoojenga_MordellWeil}. We give all details for completeness, expressing $\MW(\pi^\circ)$ in terms of sheaf cohomology and using a spectral sequence argument.

\subsection{Cohomological Characterization of $\MW(\pi^\circ)$}
Let $\mathcal{T}(\pi^\circ)$ be the sheaf of smooth translations of $\pi^\circ$: for any open set $U\subseteq B^\circ$, let $\mathcal{T}(\pi^\circ)(U)$ be the group of diffeomorphisms of $\pi^{-1}(U)$ which act by translation on each fiber.
Let $\mathcal{E}_{B^\circ}$ denote the sheaf of smooth real-valued functions on $B^\circ$. 

\begin{proposition}\label{prop: coh characterization of MW(pi^circ)}
Let $\pi\colon M\to B$ be an elliptic surface and $\pi^\circ\colon M^\circ\to B^\circ$ the fiber bundle obtained by removing a small neighborhood of each singular fiber. Then there exists a natural isomorphism
\[
    \delta\colon \MW(\pi^\circ)\xrightarrow[]{\cong} H^1(B^\circ, R^1(\pi^\circ)_*\Z),
\]
induced by the boundary map in the long exact sequence in cohomology induced by the short exact sequence of sheaves
\begin{equation} \label{eq: SES of sheaves}
    0\to R^1(\pi^\circ)_*\Z \to R^1(\pi^\circ)_*\Z\otimes_\Z \mathcal{E}_{B^\circ} \to \mathcal{T}(\pi^\circ) \to 0.
\end{equation}
\end{proposition}
\begin{proof}
For each $t\in B^\circ$, translations of $M_t$ are naturally identified with $H_1(M_t;\R)/H_1(M_t;\Z)$. Hence, fiber-wise Poincaré duality gives the short exact sequence of sheaves \eqref{eq: SES of sheaves}.
Since $\mathcal{E}_{B^\circ}$ is a fine sheaf, 
$R^1(\pi^\circ)_*\Z\otimes \mathcal{E}_{B^\circ}$ is also fine, so $H^i(B^\circ,R^1(\pi^\circ)_*\Z\otimes_\Z \mathcal{E}_{B^\circ})\cong 0$, for any $i>0$.
Hence, the long exact sequence of topological groups associated with (\ref{eq: SES of sheaves}) becomes
\[
    0\to H^0(B^\circ,R^1(\pi^\circ)_*\Z) \to H^0(B^\circ,R^1(\pi^\circ)_*\Z\otimes_\Z \mathcal{E}_{B^\circ}) \to 
    H^0(B^\circ,\mathcal{T}(\pi^\circ)) \xrightarrow[]{\delta} H^1(B^\circ,R^1(\pi^\circ)_*\Z)\to 0. 
\]
The space $H^0(B^\circ,R^1(\pi^\circ)_*\Z\otimes_\Z \mathcal{E}_{B^\circ})$ is connected, since it is a vector space, and $H^1(B^\circ,R^1(\pi^\circ)_*\Z)$ is discrete, so exactness of the above sequence implies
\[
    \pi_0(H^0(B^\circ,\mathcal{T}(\pi^\circ)))\cong H^1(B^\circ,R^1(\pi^\circ)_*\Z),
\]
which is what we wanted to prove.
\end{proof}

\subsection{$\MW(\pi^\circ)$ as annihilator of the fiber class}
For any 4-manifold $M$ and class $v\in H_2(M)$,  we denote the annihilator of $v$ by $H^2(M)^v$. We adopt the convention that a diffeomorphism $f\in\diff(M)$ acts on $H^2(M)$ by $(f^{-1})^*$.

\begin{proposition} \label{prop: MW bundle no singularities}
Let $\pi\colon M\to B$ be an elliptic fibration over a compact base $B$ and let $\pi^\circ\colon M^\circ\to B^\circ$ be the fiber bundle obtained by removing an open neighborhood of each singular fiber.
Let $F\in H_2(M^\circ)$ be the fiber class of $\pi^\circ$.
\begin{enumerate}[$(i)$]
    \item There exists a natural isomorphism
\[
    \Phi\colon\MW(\pi^\circ)  \xrightarrow[]{\cong}  H^2(M^\circ)^{F}.
\]

    \item For each $\tau\in \MW(\pi^\circ)$, there exists a unique $\gamma_\tau\in H^2(M^\circ)^{F}$ such that
$\tau$ acts on $\xi \in H^2(M^\circ)$ by
\[
    \tau(\xi)=\xi+\langle \xi,F\rangle \gamma_\tau.
\]

    \item Explicitly, the isomorphism $\Phi$ lifts to the following homomorphism 
\[
\begin{tikzcd}
    \Gamma^{\text{sm}} (B^\circ,M^\circ)\arrow[r] & H^2(M^\circ)^{F}\\
    {\sigma}\arrow[r,mapsto] &  \PD_{(M^\circ,\partial M^\circ)}[\sigma]-\PD_{(M^\circ,\partial M^\circ)}[\sigma_0|_{B^\circ}].
\end{tikzcd}
\]

\end{enumerate}

\end{proposition}
\begin{proof}
    The Leray spectral sequence for $\pi^\circ$ is
    \[
        E_2^{p,q}=H^p(B^\circ,R^q(\pi^\circ)_*\Z) \implies H^{p+q}(M^\circ).
    \]
    Since $B^\circ$ strong deformation retracts to a wedge of circles, $H^p(B^\circ,R^q(\pi^\circ)_*\Z)$ vanishes if $p\ge 2$ or $q\ge 3$. Thus, the Leray spectral sequence degenerates on the $E_2$ page.
    
    Observe that $R^2(\pi^\circ)_*\Z$ is a trivial local system of rank $1$. Indeed, every fiber of $\pi^\circ$ is a torus, so the stalks are all isomorphic to $H^2(T^2)\cong\Z$. Because the complex structure induces an orientation, $R^2(\pi^\circ)_*\Z$ must have trivial monodromy. Thus, a natural isomorphism $H^0(B^\circ,R^2(\pi^\circ)_*\Z)\cong \Z$ is given by integration over the fiber class $F\in H_2(M^\circ)$.

    Consider the Leray filtration 
    \[
        H^2(M^\circ)=L_2H^2\supset L_1H^2 \supset L_0H^2= 0.
    \]
    Since $E_2=E_\infty$, it satisfies $L_2H^2/L_1H^2 \cong H^0(B^\circ, R^2(\pi^\circ)_*\Z)$ and $L_1H^2/L_0H^2\cong H^1(B^\circ, R^1(\pi^\circ)_*\Z)$.
    This gives the short exact sequence
    \begin{equation} \label{eq: SES MW pi circ}
        0\to H^1(B^\circ, R^1(\pi^\circ)_*\Z)\to H^2(M^\circ) \to H^0(B^\circ,R^2(\pi^\circ)_*\Z) \to 0,
    \end{equation}
    where $H^0(B^\circ,R^2(\pi^\circ)_*\Z)\cong \Z$ and the second map is given by integrating over the fiber class $F\in H_2(M^\circ)$.
    The last factor is free, so the short exact sequence splits and by \Cref{prop: coh characterization of MW(pi^circ)}, 
    \[
        \MW(\pi^\circ)\cong H^1(B^\circ, R^1(\pi^\circ)_*\Z)\cong H^2(M^\circ)^{F}.
    \]

    Note that any $\tau\in \MW(\pi^\circ)$ acts on $H^2(M^\circ)$ and induces a trivial action on $H^1(B^\circ, R^1(\pi^\circ)_*\Z)$ and $H^0(B^\circ,R^2(\pi^\circ)_*\Z)$, since it acts by fiberwise translation. Let $\xi\in H^2(M^\circ)$. 
    The restriction of the zero-section $\sigma_0$ to $B^\circ$ gives a class $[\sigma_0|_{B^\circ}]\in H_2(M^\circ,\partial M^\circ)$.  
    Let $s\in H^2(M^\circ)$ denote its image under Poincaré--Lefschetz duality in $(M^\circ,\partial M^\circ)$.
    Then, $\langle s,F \rangle=1$.
    Hence, $\xi-\langle \xi, F\rangle s$ pairs trivially with $F$, so it lies in the first term of the short exact sequence (\ref{eq: SES MW pi circ}). Therefore, 
    \[
    \tau_*(\xi) = \tau_*(\xi-\langle \xi, F \rangle s) +\tau_*(\langle \xi, F\rangle s)=
    \xi-\langle \xi, F \rangle s + 
    \langle \xi, F \rangle\tau(s).
    \]
    Set $\gamma_\tau\coloneqq\tau_*(s)-s$.
    If $\tau\in \MW(\pi^\circ)$ is represented by a diffeomorphism that acts by translation by a section $\sigma$, then $\tau(s)=\PD_{(M^\circ,\partial M^\circ)}[\sigma]$.

    It remains to prove that the isomorphism $\Phi'$ defined as the composition
    \[
        \MW(\pi^\circ)\xrightarrow[]{\delta} H^1(B^\circ,R^1(\pi^\circ)_*\Z) \xrightarrow[]{\iota} H^2(M^\circ)^{F}
    \]
    is (up to sign) given by 
    \[
        \tau \mapsto \gamma_\tau.
    \]
    Let $\{U_i\}$ be a good cover of $B^\circ$ which trivializes $\pi^\circ$.
    Consider a representative $\hat\tau\in H^0(B^\circ, \mathcal{T}(\pi^\circ))$ of $\tau$ corresponding to the section $\sigma\colon B^\circ\to M^\circ$, and let $(\tau_i)_i$ be a \v{C}ech 0-cocycle representing $\hat\tau$, so $\tau_i\in \mathcal{T}(\pi^\circ)(U_i)$. Let $\lambda_i\in R^1(\pi^\circ)_*\Z\otimes_\Z \mathcal{E}_{B^\circ}(U_i)$ be lifts of $\tau_i$. On intersections, $\lambda_j-\lambda_i\in R^1(\pi^\circ)_*\Z(U_i\cap U_j)$ and
    $\delta(\hat\tau)\in H^1(B^\circ, R^1(\pi^\circ)_*\Z)$ is represented by the \v{C}ech 1-cocycle $(\lambda_j-\lambda_i)_{i,j}$. 

    We construct $b_i\in C^1(\pi^{-1}(U_i);\Z)$, such that for $i<j$, the equality $[(\lambda_i-\lambda_j)_{i,j}]=[(b_i-b_j)_{i,j}]$ holds in $H^1(B^\circ,R^1(\pi^\circ)_*\Z)$. Then the inclusion $\iota\colon H^1(B^\circ,R^1(\pi^\circ)_*\Z)\xhookrightarrow{} H^2(M^\circ)$ maps $(b_j-b_i)_{i,j}$ to the class of $(-db_i)_i$ in $H^2(M^\circ)$. Let
    \[
        q\colon H^1(T^2;\R)\to H^1(T^2;\R)/H^1(T^2;\Z)
    \]
    denote the natural projection. For each $i$, the map $\lambda_i\colon U_i\to H^1(T^2;\R)$ determines a local isotopy $p_i$ from $\sigma_0$ to $\sigma$ by setting
    \[
        p_i\colon U_i\times [0,1]\to \pi^{-1}(U_i),
    \]
    \[
        p_i(x,v)=q(v\lambda_i(x))  (\sigma_0(x)),
    \]
    where $H^1(T^2;\R)/H^1(T^2,\Z)$ is identified with the group of translations on $T^2$.
    Let $C_i\coloneqq p_i( U_i\times [0,1])$, viewed as a relative locally finite chain in $C_3^\mathrm{lf}(\pi^{-1}(U_i),\partial \pi^{-1}(U_i))$. Then $\partial C_i=[\sigma(U_i)]-[\sigma_0 (U_i)]$.
    Let $\mathcal{D}_i\colon C_k^\mathrm{lf}(\pi^{-1}(U_i),\partial \pi^{-1}(U_i))\to C^{4-k}(\pi^{-1}(U_i))$ be a chain-level realization of Poincaré--Lefschetz duality such that $d\mathcal{D}_i(c)=\mathcal{D}_i(\partial c)$ for any $c\in C_k(\pi^{-1}(U_i),\partial \pi^{-1}(U_i))$.
    Set $b_i\coloneqq \mathcal{D}_i (C_i)\in C^1(\pi^{-1}(U_i);\Z)$. Then 
    $db_i=\mathcal{D}_i([\sigma(U_i)]-[\sigma_0(U_i)])$, so this is a cocycle representing $\PD(\sigma(U_i))-\PD(\sigma_0(U_i))$.

    Note that $C_j-C_i$ is the fiberwise concatenation of the paths $p_j$ and the inverse of $p_i$. For each $x\in U_i\cap U_j$, let $l_x$ be the concatenation of the paths $p_j(x,\cdot)$ and $p_i(x,\cdot)$ with opposite orientation. This is a cycle in $M_x$, based at $\sigma_0(x)$, so by covering space theory, there exists a unique lift $\tilde l\colon [0,1]\to H^1(M_x;\R)$ with $\tilde l(0)=0$. Continuity of $\tilde l$ and the construction of $p_i$ imply 
    \[
        \tilde l(v)=\begin{cases}
            2v\lambda_j(x) \text{, if } 0\le v\le \frac{1}{2},\\
            (1-2v)\lambda_i(x)+\lambda_j(x) \text{, if } \frac{1}{2}\le v \le 1.
        \end{cases}
    \]
    Hence, $\tilde l(1)=\lambda_j-\lambda_i$ and
    $\iota[(\lambda_j-\lambda_i)_{i,j}]=\iota[(b_j-b_i)_{i,j}]=[(-db_i)_i]$, which coincide on intersections and glue to the class $-(\PD[\sigma]-\PD[\sigma_0])\in H^2(M^\circ;\Z)$.
    Hence, setting $\Phi=-\Phi'$ concludes the proof.
\end{proof}

\section{Elliptic Fibrations over a Disk with One Singular Fiber} \label{section: local picture}
In this section, we analyze an elliptic fibration over a disk with one singular fiber, explicitly computing $\MW(\pi)$ and $\MW(\pi^\circ)$.
We determine a complete obstruction to extending sections locally, in a small neighborhood of a singular fiber. 
Let $D$ be the closed unit disk in $\C$ and consider an elliptic fibration $\pi\colon X\to D$ with only one singular fiber at the origin. Let $D_\varepsilon\subset D$ be a small open disk around the origin and let $D^\circ\coloneqq D-D_\varepsilon$ and $X^\circ\coloneqq \pi^{-1}(D^\circ)$.  We answer the following: 
\begin{question}
    Given an isotopy class in $\MW(\pi^\circ)$, when does it have a representative section $\sigma\colon D^\circ\to X^\circ$ that extends to a section $\sigma\colon D\to X$?
\end{question}
Since $D^\circ$ strong deformation retracts onto its boundary, this is equivalent to asking when a section 
$\sigma\colon \partial D\to \partial X$ extends to a section over the entire disk $D$. We construct a map $\psi\colon \MW(\pi^\circ)\to H_1(X_0)$ and show that it captures the obstruction to extending sections: a section $\sigma$ extends if and only if $\psi([\sigma])=0$.

\subsection{Definition of $\psi$ and Local Obstruction Theorem}
\begin{figure}[h]
    \centering
    \fbox{\includegraphics[width=0.31\textwidth]{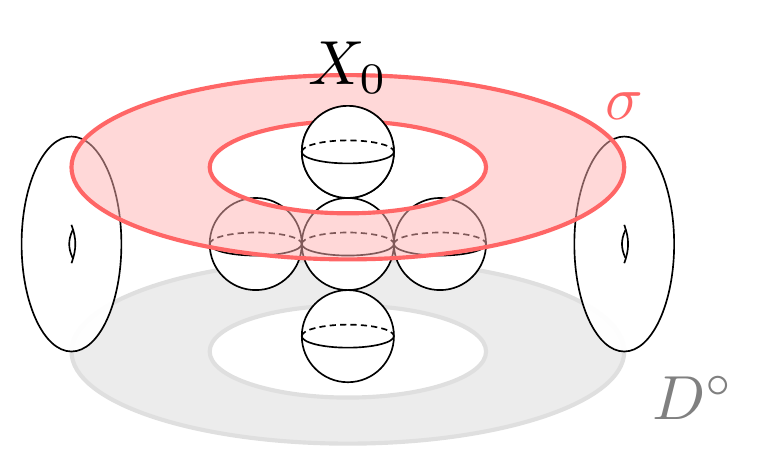}}\hfill%
    \fbox{\includegraphics[width=0.31\textwidth]{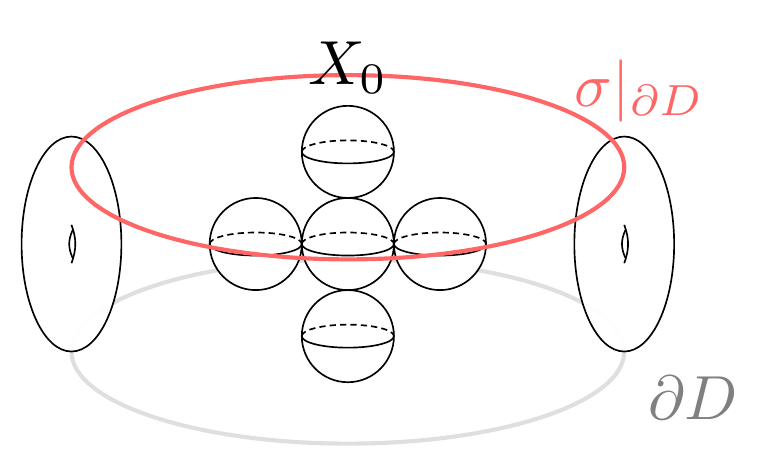}}\hfill%
    \fbox{\includegraphics[width=0.31\textwidth]{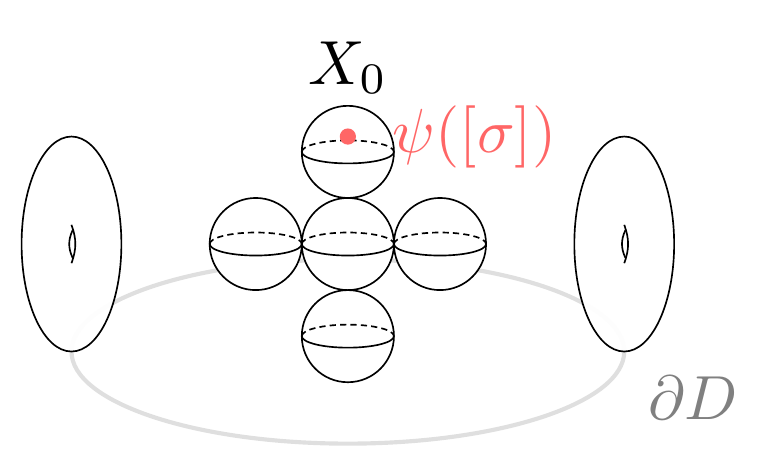}}
    
    \caption{The map $\psi$ restricts a section $\sigma$ to the boundary, then shrinks it to $X_0$.}
    \label{fig:section to H1}
\end{figure}

We define a homomorphism $\psi\colon \MW(\pi^\circ)\to H_1(X_0)$. Given an element in $\MW(\pi^\circ)$, pick a representative  section $\sigma\colon D^\circ\to X^\circ$. The image of the restriction $\sigma|_{\partial D}\colon \partial D\to \partial X$ represents a homology class ${[\sigma|_{\partial D}]\in H_1(\partial X)}$. This class is independent of the choice of isotopy representative. Indeed, if two sections $\sigma,\sigma'\colon D^\circ\to X^\circ$ are isotopic, the isotopy restricts to $\partial D$ and traces an annulus in $\partial X$. Then, $\sigma(\partial D)$ and $\sigma'(\partial D)$ are the two boundary components of the annulus traced out by the isotopy. Thus, the 1-cycle in $\partial X$ represented by $\sigma|_{\partial D}-\sigma'|_{\partial D}$ is the boundary of a $2$-cycle in $\partial X$, meaning $[\sigma|_{\partial D}]=[\sigma'|_{\partial D}]$ in $H_1(\partial X)$. Hence, we can define a map  
\[
    \eta\colon \MW(\pi^\circ)\to H_1(\partial X),
\]
\[
    [\sigma]\mapsto [\sigma|_{\partial D}]-[\sigma_0|_{\partial D}].
\]
The inclusion $\partial X\xhookrightarrow{}X$ induces a map
$H_1(\partial X) \to H_1(X)$. By \Cref{prop: noighborhood of singular fiber strong deformation retracts}, $X$ strong deformation retracts to $X_0$, inducing an isomorphism $H_1(X)\cong H_1(X_0)$. We define $\psi$ as the composition of these maps:
\begin{center}
\begin{tikzcd}
    \MW(\pi^\circ) \arrow[r,"\eta"] \arrow[rrr, bend right=20,"\psi"]  & 
    H_1(\partial X) \arrow[r] & H_1(X) \arrow[r,"\cong"] & H_1(X_0).
\end{tikzcd}
\end{center}
Intuitively, $\psi$ takes a section over $D^\circ$ and pushes its image over $\partial D$ into the central fiber $X_0$, as depicted in Figure \ref{fig:section to H1}.

The following theorem gives a complete characterization of when a section extends to the central fiber up to isotopy. This is crucial in computing the Mordell--Weil group of an arbitrary elliptic surface.

\begin{theorem}[\textbf{Local Obstruction Theorem}] \label{thm: extend sections locally}
    Let $\pi\colon X\to D$ be an elliptic fibration over a disk, with only one singular fiber at the origin. 
    The following is a short exact sequence of groups
    \[
        0\to \MW(\pi)\to \MW(\pi^\circ)\xrightarrow[]{\psi} H_1(X_0) \to 0.
    \]
\end{theorem}

\begin{remark} \label{rem: extend implies zero in homology of center}
    We can see immediately that if a section $\sigma\colon D^\circ\to X^\circ$ extends to a section $\Bar{\sigma}\colon D\to X$, then $\psi([\sigma])=0$. Indeed, $\sigma|_{\partial D}$ represents a trivial class in $H_1(X)$, since it is the boundary of the $2$-cycle $\Bar{\sigma}(D)$. The same holds for the zero-section $\sigma_0$. Therefore, $[\sigma|_{\partial D}]-[\sigma_0|_{\partial D}]$ is zero in $H_1(X)\cong H_1(X_0)$.
    The main point of Theorem \ref{thm: extend sections locally} is that this condition is sufficient: if $\psi([\sigma])=0$, then $[\sigma]$ has a representative which extends to $D$.
\end{remark}

\begin{remark}
    Recall that by Kodaira's classification, additive fibers are simply-connected. Hence, if $X_0$ is additive, $H_1(X_0)\cong 0$ and Theorem \ref{thm: extend sections locally} implies that all section classes of $\pi^\circ$ extend to sections of $\pi$.
\end{remark}

\begin{figure}
    \centering
    \begin{subfigure}[b]{0.49\textwidth}
        \centering
        \includegraphics[width=\textwidth]{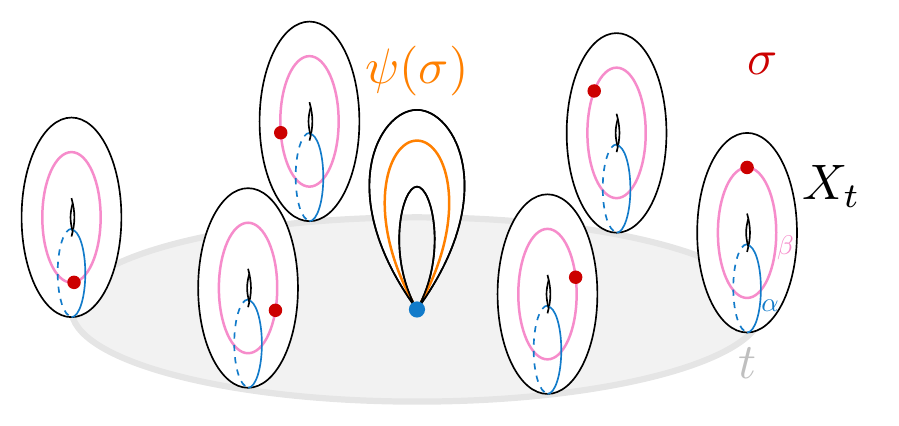}
        \caption{$\psi$ maps the isotopy class of the section $\sigma\colon \partial D \to \partial X$ constructed in Example \ref{ex: section does not extend} to $\psi(\sigma)$ in $H_1(X_0)$. This is a non-trivial homology class, so $\sigma$ does not extend to $D$.}
        \label{fig:ex section does not extend}
    \end{subfigure}
    \hfill
    \begin{subfigure}[b]{0.49\textwidth}
        \centering
        \includegraphics[width=\textwidth]{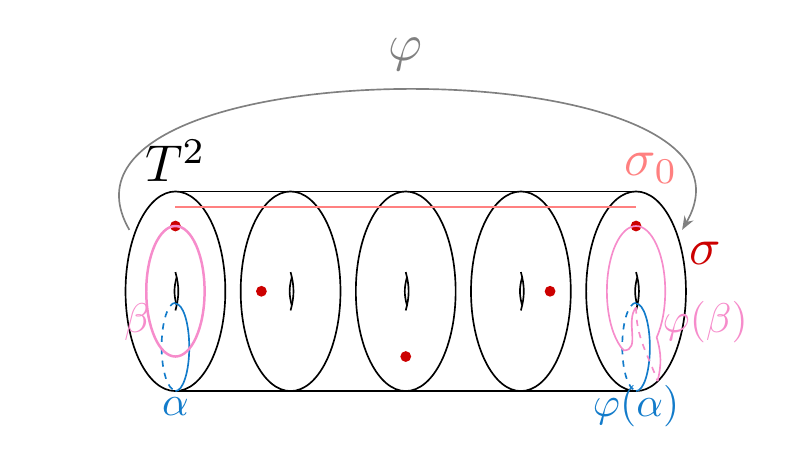}
        \caption{$\partial X$ is diffeomorphic to the mapping torus $M_\varphi$.}
        \label{fig:mapping torus}
    \end{subfigure}
    \caption{Example of a section that does not extend.}
    \label{fig: ex section not extend both images}
\end{figure}

\begin{ex}[\textbf{A section that does not extend}]\label{ex: section does not extend}
    Let $D\subset \C$ be a disk centered at the origin and let $\pi\colon X\to D$ be a proper holomorphic map such that for $t\ne 0$, the fibers $X_t$ are smooth tori and the fiber at the origin is of type $I_n$\footnote{This is the same as an elliptic fibration without a fixed zero-section.}  (see Figure \ref{fig:ex section does not extend} for the case $n=1$). The restriction $\pi\colon \partial X\to \partial D$ is a torus-bundle over a circle. Pick a basepoint $t$ on $\partial D$ and fix an identification of $X_t$ with the standard torus $T^2$.
    Choose a parametrization $\gamma\colon [0,1]\to \partial D$ with $\gamma(0)=\gamma(1)=t$.  Let $f\in \Mod(T^2)$ be the monodromy transformation corresponding to $\gamma$. By Kodaira's classification, we can pick curves $\alpha,\beta$ on $T^2$, representing a basis of $H_1(T^2)$ such that $f$ acts on $H_1(T^2)$ as
    \[ 
        f_*=\begin{pmatrix}
            1 & n\\
            0&1
        \end{pmatrix}.
    \]
    This means that $f$ is the $n$-th power of a Dehn twist about $\alpha$. By Picard-Lefschetz theory, this implies that $\alpha$ is the vanishing cycle of the fibration, i.e.\ its homology class vanishes in $H_1(X_0)$. Let $\varphi\in \mathrm{Diff}^+(T^2)$ be a representative of $f$ that acts as the identity on $T^2-R$ for a regular neighborhood $R$ of $\alpha$.
    Then $\partial X$ is diffeomorphic to the mapping torus 
    \[
        M_\varphi=(T^2\times [0,1])/(x,0)\sim (\varphi(x),1),
    \]
    see \Cref{fig:mapping torus}.
    
    In order to obtain an elliptic fibration,
    we fix a zero-section $\sigma_0\colon D\to X$. Using the mapping torus description of $\partial X$, we define $\sigma_0$ over $\partial D$, then find an explicit extension. Pick $x\in T^2-R$. Define $\sigma_0'\colon [0,1]\to M_\varphi$ by $\sigma_0'(s)=[(x,s)]$. Since $x\in T^2-R$, it follows $\varphi(x)=x$ and thus $\sigma_0'(0)=\sigma_0'(1)$. Hence, $\sigma_0'$ induces a map $\sigma_0'\colon \partial D\to M_\varphi$. Composing with the diffeomorphism $M_\varphi\to \partial X$ gives a section $\sigma_0\colon \partial D\to \partial X$. This section picks out ``the same'' point in each fiber, which is away from the vanishing cycle. Thus, we can extend it radially to $D$. Since $\sigma_0|_{\partial D}$ never meets the vanishing cycle, $\sigma_0(0)$ will not meet the node of $X_0$. Hence, $\sigma_0(0)\in X_0^\#$. Thus, $\sigma_0\colon D\to X$ is a section of $\pi$.

    We  define a section $\sigma\colon \partial D\to \partial X$ and show it does not extend to $D$ (see Figure \ref{fig: ex section not extend both images}). Choose a parametrization $\beta\colon [0,1]\to T^2$ such that $\beta(0)=\beta(1)\in T^2-R$.
    Define $\sigma\colon [0,1]\to M_\varphi$, $\sigma(s)=[(\beta(s),s)]$. Since $\beta(0)\in T^2-R$, it is fixed by $\varphi$. Thus, 
    \[
        \sigma(0)=[(\beta(0),0)]=[(\varphi(\beta(0)),1)]= [(\beta(0),1)]= [(\beta(1),1)]=\sigma(1).
    \]
    Hence, this induces a map $\sigma\colon \partial D\to M_\varphi$.
    Composing it with the diffeomorphism $M_\varphi\to\partial X$ defines a section of $\pi\colon \partial X\to \partial D$.

    The map $\psi$ takes $[\sigma]$ to the loop $[\beta]\in H_1(X_0)$. Recall that $\alpha$ is the vanishing cycle, so $[\beta]$ is the generator of $H_1(X_0)\cong \Z$.
    By Remark \ref{rem: extend implies zero in homology of center}, $\sigma$ cannot extend over the whole disk, since $[\beta]$ is non-trivial in $H_1(X_0)$.
\end{ex}

\subsection{Computation of $\MW(\pi)$}
In the setup of an elliptic fibration over a disk with only one singular fiber, choosing the isotopy class of a section is equivalent to choosing a multiplicity one component of the singular fiber. We make this more explicit in the following proposition.
\begin{proposition}[\textbf{Computation of $\MW(\pi)$}] \label{prop: MW(pi) if only one singularity}
Let $\pi\colon X\to D$ be an elliptic fibration over a disk with only one singular fiber at $0$. Let $X_0^\#$ be the smooth locus of $X_0$ and let $G(X_0)\coloneqq \pi_0(X_0^\#)$ be the group of components of the regular part of the singular fiber $X_0$
\footnote{Recall that this means that $G(X_0)$ is the group of irreducible components of $X_0$ of multiplicity $1$. The group structure on $X_0^\#$, is $\C \times G(X_0)$ or $\C^\times \times G(X_0)$.}. Denote by $\Gamma^{\mathrm{sm}}(\pi)$ the group of smooth sections of $\pi$, where we choose $\sigma_0$ to be the identity element.
There is a homomorphism $\varphi\colon \Gamma^{\mathrm{sm}}(\pi)\to G(X_0)$
defined as the composition 
\[
\begin{tikzcd}
        \Gamma^{\mathrm{sm}}(\pi)\arrow[r]\arrow[rr, bend right=20, "\varphi"'] & X_0^\# \arrow[r] & G(X_0),
\end{tikzcd}
\]
where the first map evaluates a section $\sigma\colon D\to X$ at $0$ and the second is projection $X_0^\#\to \pi_0(X_0^\#)$.

Furthermore, $\varphi$ induces an isomorphism 
\[
        \Bar{\varphi}\colon \MW(\pi)\xrightarrow[]{\cong} G(X_0).
    \]
\end{proposition}
\begin{proof}
    \textit{$\varphi$ is well-defined.}
    Consider a section $\sigma\colon D\to X$. By definition, $\pi\circ\sigma=\id_D$, so the chain rule implies that $d\pi_0$ has full rank. Thus, $\sigma(0)\in X_0^\#$, so $\varphi$ is well-defined. Furthermore, $\varphi$ is a group homomorphism, as a composition of group homomorphisms.

    \textit{$\varphi$ does not depend on the isotopy representative.}
    Given two isotopic sections $\sigma,\sigma^1\colon D \to X$, there exists a one-parameter family of sections $\sigma^t\colon D\to X$ with $\sigma^0=\sigma$, for $t\in [0,1]$. As we showed above, $\sigma^t(0)\in X_0^\#$. By continuity in $t$, it follows that $\sigma^t(0)$ lies in the same connected component of $X_0^\#$ for all $t$. In particular, $\sigma$ and $\sigma^1$ intersect the same irreducible component of $X_0$, inducing a homomorphism $\MW(\pi)\to G(X_0)$.

    \textit{$\varphi$ is surjective}. Let $C_i\in G(X_0)$ be a component of $X_0$ of multiplicity $1$ and pick a regular point $x\in C_i$.
    Since $\pi$ is a submersion at $x$ there exist local coordinates around $x\in X$ and $0\in D$ such that $\pi$ is the projection $(z,w)\mapsto z$. 
    In these coordinates, define a section $\sigma(z)=(z,0)$. Let $D_\varepsilon\subset D$ be an open disk centered at the origin, completely contained in the coordinate chart on $D$ chosen above.
    Let $X_\varepsilon\coloneqq\pi^{-1}(D_\varepsilon)$. Then $\pi\colon X-X_\varepsilon\to D-D_\varepsilon$ is a torus-bundle and $\sigma|_{\partial D_\varepsilon}$ is a well-defined section over $\partial D_\varepsilon$. Since the annulus $D-D_\varepsilon$ strong deformation retracts to $\partial D_\varepsilon$,
    the bundle $X-X_\varepsilon\to D-D_\varepsilon$ is isomorphic to the pullback of $\pi^{-1}(\partial D_\varepsilon)\to \partial D_\varepsilon$ along the retraction. Therefore, $\sigma|_{\partial D_\varepsilon}$ can be extended to a section over the entire annulus, yielding a section $\sigma\colon D\to X$.
    
    \textit{$\Bar{\varphi}$ is injective}. Let $\sigma\colon D\to X$ be a section intersecting the zero-component $C_0$ of $X_0$. 
    We first use a flow argument to isotope $\sigma$ to a section that meets $X_0$ at $\sigma_0(0)$. Then we use the fact that $\pi$ is a submersion on the smooth locus of $X$ to conclude.
    
    Let $C_0^\#$ denote the zero-component with singular points removed. By Kodaira's classification, the fact that $X_0$ is singular implies that $C_0^\#$ is diffeomorphic to $\C$ or $\C^\times$. In particular, it is path-connected, so there exists a path $\gamma\colon [0,1]\to C_0^\#$ with $\gamma(0)=\sigma(0)$ and $\gamma(1)=\sigma_0(0)$. Assume that $\gamma$ is an embedding and choose a compactly supported vector field $v$ on $C_0^\#$ such that $\gamma'(t)=v_{\gamma(t)}$, for all $t\in [0,1]$. The vector field $v$ can be smoothly extended to a compactly supported vector field $V$ on an open neighborhood $W\subset X^\#$. We extend $V$ to the entire $X^\#$ by setting it equal to $0$ everywhere outside $W$. Let $T_VX^\#\coloneqq \ker(d\pi)$ be the vertical tangent bundle. 
    Fix a Riemannian metric on the manifold $X^\#$ and let
    $\tilde{V}_x$ denote the orthogonal 
    projection of $V_x$ to $T_VX^\#_x$.
    Recall that we defined $v$ as a vector field on $C_0^\#\subset X_0^\#$, so $d\pi(v)=0$, meaning $v\in T_VX^\#$. Thus, $\tilde V|_{C_0^\#}=v$.

    Integrating the compactly supported vector field $\tilde V$ gives a global flow $\Phi_t\colon X^\#\to X^\#$. 
    For arbitrary $t\in[0,1]$, we have
    \[
        \frac{d}{dt}(\pi\circ \Phi_t(x)) = d\pi_{\Phi_t(x)}\left(\frac{d}{dt}\Phi_t(x)\right) = d\pi_{\Phi_t(x)}(\tilde V_{\Phi_t(x)})=0.
    \]
    Thus, $\pi\circ \Phi_t(x)$ is constant in $t$, meaning $\pi\circ \Phi_t(x)=\pi\circ \Phi_0(x)=\pi(x)$. Hence, $\Phi_t$ maps each fiber to itself. In particular, it maps sections to sections. Let $\tilde\sigma\colon D\to X$ be defined by $\tilde \sigma\coloneqq \Phi_1\circ\sigma$. This is a section of $\pi$ with $\tilde\sigma(0)=\Phi_1(\sigma(0))=\sigma_0(0)$, by definition of the vector field $v$. Define an isotopy $H\colon D\times [0,1]\to X$ by setting $H(s,t)=\Phi_t(\sigma(s))$. Thus, $\sigma$ and $\tilde \sigma$ are isotopic sections.

    It remains to show that $\tilde \sigma$ and $\sigma_0$ are isotopic. 
    Let $p\coloneqq \sigma_0(0)=\tilde\sigma(0)$. This is a regular point, so by the local submersion theorem, there exists a chart $U$ on $X$ around $p$ and $D_\varepsilon$ on $D$ around 0, such that in this chart, $\pi$ is just a projection. We can define an isotopy over $D_\varepsilon$ using a straight-line homotopy in the local chart. Now observe that the obstruction to extending an isotopy over the annulus $D-D_\varepsilon$ lies in $H^k(D-D_\varepsilon,\partial D_\varepsilon;\pi_k(T^2))$,  which is trivial since the annulus $D-D_\varepsilon$ strong deformation retracts to $\partial D_\varepsilon$. By \cite[Theorem 34.6]{Steenrod_FiberBundles}, the isotopy over $\partial D_\varepsilon$ extends over the entire annulus $D-D_\varepsilon$.
    \end{proof}

    \begin{corollary}\label{cor: MW pi is finite}
        If $\pi\colon X\to D$ is an elliptic fibration with only one singular fiber, then $\MW(\pi)$ is finite.
    \end{corollary}
    \begin{proof}
        By Proposition \ref{prop: MW(pi) if only one singularity}, $\MW(\pi)\cong G(X_0)$, which is finite by Kodaira's classification.
    \end{proof}

\subsection{Explicit Computation of $\MW(\pi^\circ)$}
We give a homological characterization of $\MW(\pi^\circ)$, identifying it with the group of monodromy coinvariants.
Fix a basepoint $t\in \partial D$.
Let $\gamma$ be a counterclockwise generator for $\pi_1(\partial D,t)$. Let $f\in \Mod(X_t)$ be the monodromy transformation corresponding to the loop $\gamma$. 
We show that the homomorphism 
\[
    \eta\colon \MW(\pi^\circ)\to H_1(\partial X),
\]
\[
    [\sigma]\mapsto [\sigma|_{\partial D}] - [\sigma_0|_{\partial D}]
\]
is injective, identifying $\MW(\pi^\circ)$ with a specific subgroup of $H_1(\partial X)$.
The following lemma also shows that $H_1(\partial X)$ consists of ``vertical elements'' (elements in $\MW(\pi^\circ)$) and ``horizontal elements'' (elements of $H_1(S^1)$).

\begin{lemma}[\textbf{Computation of $\MW(\pi^\circ)$}] \label{lem: sections bij to coinvariants monodromy}
    Let $\pi\colon X\to D$ be an elliptic fibration over a disk with only one singular fiber at $0$ and monodromy $f\in \Mod(X_t)$. The Mordell--Weil group of $\pi^\circ$ fits into the following short exact sequence
    \begin{equation}\label{eq: SES MW pi^circ and H_1partial X}
        0\to \MW(\pi^\circ) \xrightarrow[]{\eta} H_1(\partial X) \xrightarrow[]{\pi_*} H_1(\partial D) \to 0.
    \end{equation}

    Furthermore,
    \[
        \MW(\pi^\circ)\cong \frac{H_1(X_t)}{(I-f_*)H_1(X_t)} \cong \Z^2/(I-f_*)\Z^2.
    \]
\end{lemma}
\begin{proof}
    Let $F\coloneqq [X_t]\in H_2(X^\circ)$ denote the fiber class of $\pi^\circ$ and let $H^2(X^\circ)^F$ denote the annihilator of $F$.
    By Proposition \ref{prop: MW bundle no singularities}, there exists an isomorphism $\MW(\pi^\circ)\to H^2(X^\circ)^F$, given by $[\sigma]\mapsto \PD[\sigma]-\PD[\sigma_0]$.
    Since $X^\circ$ strong deformation retracts onto $\partial X$, restriction gives an isomorphism $H^2(X^\circ)^F\xrightarrow{} H^2(\partial X)^F$. Poincaré duality on the closed oriented three-manifold $\partial X$ gives an isomorphism $\PD_{\partial X}\colon H^2(\partial X)\xrightarrow[]{} H_1(\partial X)$. Thus, composing these maps gives an inclusion 
    \[
        \MW(\pi^\circ)\xhookrightarrow{} H_1(\partial X),
    \]
    \[
        [\sigma]\mapsto [\sigma|_{\partial D}]-[\sigma_0|_{\partial D}].
    \]
    The image of this inclusion is equal to the image of $H^2(\partial X)^F$ under Poincaré duality.
    A cohomology class $\alpha\in H^2(\partial X)$ lies in the annihilator of the fiber class $F$ if and only if $0=\langle \alpha,F\rangle = \PD_{\partial X}(\alpha)\cdot F$, where $\cdot$ denotes the algebraic intersection number. Hence, 
    \[
        \PD_{\partial X}(H^2(\partial X)^F) = \{a\in H_1(\partial X)\mid  a\cdot F=0 \}.
    \]

    Consider the Wang exact sequence for the fiber bundle $X_t\to \partial X \to S^1$
    \begin{equation} \label{eq: Wang mapping torus}
        \dots\to H_1(X_t)\xrightarrow[]{I-f_*} H_1(X_t) \xrightarrow[]{} H_1(\partial X) \xrightarrow[]{\partial} H_0(X_t)\xrightarrow[]{I-f_*} H_0(X_t)\to \dots
    \end{equation}

    The monodromy $f$ acts as the identity on $H_0(X_t)$, so the map $H_0(X_t)\xrightarrow[]{I-f_*} H_0(X_t)$ in the long exact sequence \eqref{eq: Wang mapping torus} is zero. Thus, we obtain a short exact sequence
    \begin{equation}\label{eq: SES H_1 partial X only one singularity}
        0\to \frac{H_1(X_t)}{(I-f_*)H_1(X_t)} \to H_1(\partial X)\xrightarrow[]{\partial} H_0(X_t)\to 0.
    \end{equation}

    The boundary map $H_1(\partial X) \xrightarrow[]{\partial} H_0(X_t)$ is precisely the algebraic intersection number with the class $[X_t]\in H_2(\partial X)$. Hence, the image of $H^2(\partial X)^F$ under Poincaré duality in $\partial X$ is equal to the kernel of $H_1(\partial X) \xrightarrow[]{\partial} H_0(X_t)$.
    Therefore,
    
    \begin{align*}
        &\MW(\pi^\circ)\cong H^2(X^\circ)^F\cong \PD_{\partial X}(H^2(\partial X)^F)=\\
        &=\ker(H_1(\partial X)\xrightarrow[]{\partial} H_0(X_t))\cong  \frac{H_1(X_t)}{(I-f_*)H_1(X_t)}\cong  \Z^2/(I-f_*)\Z^2.
    \end{align*}
    
    We finish the proof by showing that the short exact sequence \eqref{eq: SES H_1 partial X only one singularity} implies exactness of the short exact sequence \eqref{eq: SES MW pi^circ and H_1partial X}. As explained before, the map $\partial\colon H_1(\partial X)\to H_0(X_t)\cong \Z$ is the algebraic intersection number with the class $[X_t]\in H_2(\partial X)$. This is equivalent to saying that $\partial$ measures how many times a homology class winds around the circle $\partial D$. Hence, under the canonical identifications $H_0(X_t)\cong \Z$ and $H_1(\partial D)\cong\Z$, we get $\forall c\in H_1(\partial X)$, $\partial c= c\cdot [X_t] = (\pi|_{\partial X})_*(c)$. 
    Hence, $\partial$ can be replaced with $(\pi|_{\partial X})_*\colon H_1(\partial X)\to H_1(\partial D)\cong \Z$.
\end{proof}

\begin{corollary}\label{cor: additive implies MW boundary finite}
    Let $\pi\colon X\to D$ be an elliptic fibration with only one singular fiber at the origin. The singular fiber $X_0$ is additive if and only if $\MW(\pi^\circ)$ is finite.
\end{corollary}
\begin{proof}
    By Lemma \ref{lem: sections bij to coinvariants monodromy}, $\MW(\pi^\circ)\cong \Z^2/(I-f_*)\Z^2$. This group is finite if and only if $\det(I-f_*)\ne 0$. Thus, it suffices to show that $1$ is an eigenvalue of $f_*\colon H_1(X_t)\to H_1(X_t)$ if and only if $X_0$ is not additive. The determinant of $f_*$ is $1$, so if $1$ is an eigenvalue of $f_*$, then both eigenvalues must equal $1$. Therefore, $\mathrm{tr}(f_*)=2$. By Kodaira's classification, the trace of the monodromy is $2$ if and only if $X_0$ is smooth or multiplicative.
\end{proof}

\subsection{Injectivity of the restriction map $\MW(\pi)\to \MW(\pi^\circ)$}
We show that, if two sections of $\pi$ are isotopic over $\partial D$, then they must be isotopic over $D$.

\begin{figure}
    \centering
    \includegraphics[width=0.4\linewidth]{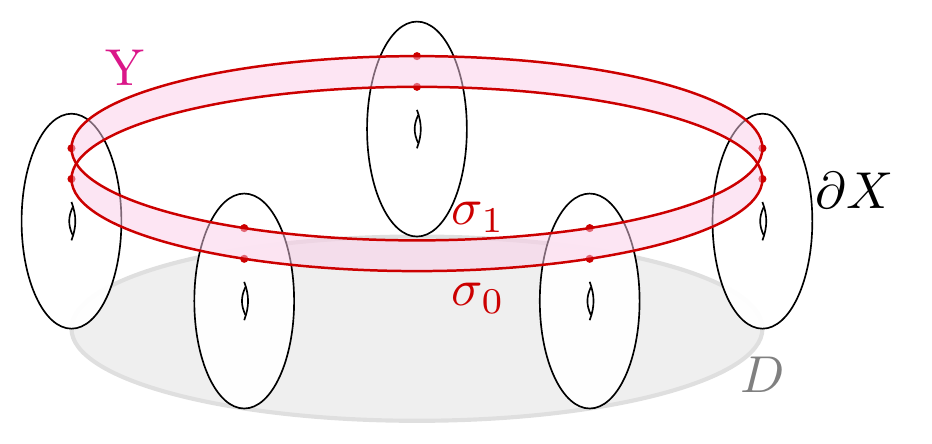}
    \caption{Two isotopic sections $\sigma,\sigma_0$ trace out an annulus $Y$ in $\partial X$.}
    \label{fig:isotopy gives annulus}
\end{figure}

\begin{proposition}\label{prop: isotopic on circle implies isotopic on disk}
    Let $\pi\colon X\to D$ be an elliptic fibration over a disk $D$ with only one singular fiber at the origin. Then the homomorphism $\MW(\pi)\to \MW(\pi^\circ)$ given by restricting sections to $D^\circ$ is injective.
\end{proposition}
\begin{proof}
    Let $\sigma\colon D\to X$ be a section of $\pi$ such that $\sigma|_{\partial D}\simeq \sigma_0|_{\partial D}$.
    By Proposition \ref{prop: MW(pi) if only one singularity}, $\sigma$ and $\sigma_0$ are isotopic over $D$ if and only if they intersect the same component of $X_0$. Let us assume that they intersect different components of $X_0$.
    Denote the irreducible components of $X_0$ by $C_0,\dots,C_r$ and assume that $\sigma_0$ intersects $C_0$ and $\sigma$ intersects $C_a$ for $a\ne 0$. Since sections intersect each fiber with multiplicity $1$, $C_a$ must be a multiplicity $1$ component of $X_0$. Let $-\sigma(D)$ be $\sigma(D)$ with reversed orientation. The algebraic intersection number satisfies $-\sigma(D)\cdot [X_0]=-1$.

    The isotopy between $\sigma|_{\partial D}$ and $\sigma_0|_{\partial D}$ traces out an annulus $Y$ in $\partial X$, see Figure \ref{fig:isotopy gives annulus}. Since $\partial \sigma_0(D)$ and $\partial(-\sigma(D))$ have opposite orientation in $\partial X$, they can be glued to the respective boundaries of $Y$, yielding 
    the closed surface $S\coloneqq \sigma_0(D)\cup (-\sigma(D)) \cup Y$. This defines a homology class $[S]\in H_2(X)$. Since $H_2(X)\cong H_2(X_0)$, which is freely generated by the irreducible components of $X_0$, there exist $x_i\in \Z$ such that $[S]=\sum x_i [C_i]$. Since $Y\subset \partial X$, it does not intersect $X_0$. This gives the following system of equations:
    \begin{equation*}
    \begin{cases}
        1&=[S]\cdot [C_0] = \sum x_i[C_i][C_0] \\
            -1&=[S]\cdot [C_a] = \sum x_i[C_i][C_a] \\
            0&=[S]\cdot [C_j] = \sum x_i[C_i][C_j], \text{ for } j\ne0,a.
    \end{cases}  
    \end{equation*}
    Let $\tilde A\in \mathrm{M}_{(r+1)\times (r+1)}(\Z)$ be the extended Cartan matrix associated with the singular fiber $X_0$.\footnote{The matrix $-\tilde{A}$ encodes pairwise intersection of components of $X_0$.} 
    The system of equations is equivalent to the matrix equation $-\tilde Ax=e_0-e_a$, where $x=(x_0,\dots,x_r)^T$ and $e_0$ and $e_a$ denote the zeroth and $a$-th standard basis vector in $\R^{r+1}$, respectively.
    Let $A$ be the Cartan matrix obtained by removing the component $C_0$. We can write
    \[\tilde A=
    \begin{pmatrix}
        2 & v^T \\
        v & A
    \end{pmatrix},
    \]
    where $v\in \Z^{r}$.
    Let $\tilde m\in \Z^{r+1}$ be the vector of multiplicities of the irreducible components of $X_0$.  By Kodaira's classification, this vector spans the kernel of $\tilde A$. Since $C_0$ has multiplicity $1$, we can write $\tilde m=\begin{pmatrix}
        1&m
    \end{pmatrix}^T$ for $m\in \Z^{r}$. Now, $\tilde A \tilde m=0$ implies $v+Am=0$. Also, writing $x=\begin{pmatrix}
        x_0 &x'
    \end{pmatrix}^T$ with $x'\in\Z^{r}$ in the equation $\tilde Ax=-e_0+e_a$ yields
    \[
        x_0v+Ax'=\bar{e}_a,
    \]
    where $\bar{e}_a$ denotes the $a$-th standard basis vector in $\Z^r$.
    Using $v+Am=0$, this implies
    \begin{equation} \label{eq: cartan matrix computation}
        A(x'-x_0m)=\bar{e}_a.
    \end{equation}
    The Cartan matrix $A$ is invertible, so \eqref{eq: cartan matrix computation} is equivalent to
    \begin{equation}\label{eq: cartan matrix with inverse}
        x'-x_0m = A^{-1}\bar{e}_a.
    \end{equation}

    In what follows, we show that the $a$-th column of the inverse of a Cartan matrix of type ADE cannot be integral.
    Since $x',m\in \Z^r$ and $x_0\in\Z$, this implies that \eqref{eq: cartan matrix with inverse} is impossible and gives the desired contradiction.
    
    Recall that $a$ corresponds to a non-zero component of multiplicity one, so we only need to check such columns. We do this explicitly, by considering the inverses of Cartan matrices of each type. We check that the columns corresponding to multiplicity $1$ components are not integral. All formulas for inverses of Cartan matrices we use below are from \cite{WeiZou_InversesCartan}.
    \begin{itemize}
        \item \textit{Type $A_n$:} All components have multiplicity $1$;
        we need to check all columns. If $A$ is a Cartan matrix of type $A_n$, then
        \[
        A_{i,j}=\begin{cases}
            2 \text{, if } i=j,\\
            -1 \text{, if } |i-j|=1,\\
            0 \text{, otherwise.}           
        \end{cases}
        \]
        Its inverse has entries 
        \[
            (A^{-1})_{i,j}= \min(i,j) -  \frac{ij}{n+1} \text{, for } 1\le i,j\le n.  
        \]
        For each $1\le j\le n$, the first entry of the $j$-th column is $1-\frac{j}{n+1}\notin \Z$.
        \item \textit{Type $D_n$:} In this case we assume $n\ge 4$. We order the vertices as in Figure \ref{fig:Dynkin with labels}. Only the first and the last two components have multiplicity $1$, so $a\in\{1,n-1,n\}$. The corresponding Cartan matrix is
        \[
            A=\begin{pmatrix}
                2 & -1 & & & &\\
                -1 &\ddots& \ddots &&&\\
                &\ddots& 2&-1&-1\\
                &&-1&2&0\\
                &&-1&0&2\\                
            \end{pmatrix}
        \]
        The matrix is symmetric and so is its inverse, so it suffices to give the entries of the inverse for $i\le j$:
        \[
            (A^{-1})_{i,j}=\begin{cases}
                i, \text{ if } 1\le i\le j\le n-2\\
                i/2, \text{ if } i<n-1, j=n-1 \text{ or } n\\
                \frac{n-2}{4}, \text{ if } i=n-1,j=n\\
                \frac{n}{4}, \text{ if } i=j=n-1 \text{ or } n.\\
            \end{cases}
        \]
        We start with the column given by $j=1$. By symmetry, we have $(A^{-1})_{n,1}=(A^{-1})_{1,n}=1/2$, so this column is not integral.
        For the columns $j=n-1$ and $j=n$, the fact that $n\ge 4$ implies that $(A^{-1})_{1,n-1}=(A^{-1})_{1,n}=1/2$, so these cannot be integral.

    \begin{figure}
    \centering
    \includegraphics[width=1\linewidth]{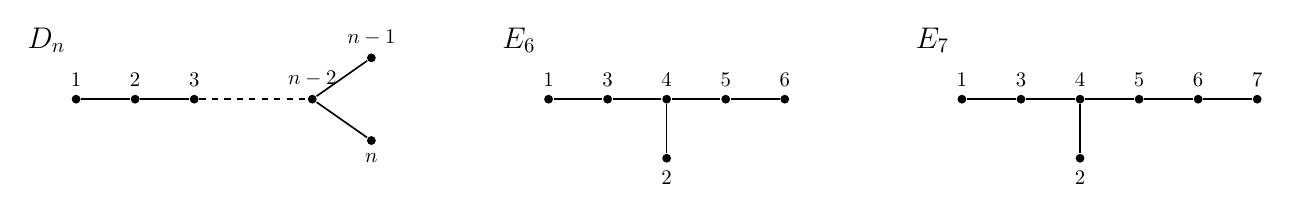}
    \caption{Labelling of vertices for Dynkin diagrams.}
    \label{fig:Dynkin with labels}
    \end{figure}

        \item \textit{Type $E$:}
        \begin{itemize}
            \item $E_6$: If we order the vertices as in Figure \ref{fig:Dynkin with labels}, then only the first and sixth component have multiplicity $1$, so $a\in\{1,6\}$. The corresponding Cartan matrix is
            \[
                A=\begin{pmatrix}
                    2&0&-1&0&0&0\\
                    0&2&0&-1&0&0\\
                    -1&0&2&-1&0&0\\
                    0&-1&-1&2&-1&0\\
                    0&0&0&-1&2&-1\\
                    0&0&0&0&-1&2
                \end{pmatrix}.
            \]
            We only need to give a non-integral entry in the first column of the inverse and one in the last column of the inverse. $(A^{-1})_{1,1}=4/3$ and $(A^{-1})_{1,6}=2/3$, so we are done.
            \item $E_7$: If we order the vertices as in Figure \ref{fig:Dynkin with labels}, only the seventh component has multiplicity $1$, so $a=7$. Hence, we need to check only the last column of the inverse of the Cartan matrix. $(A^{-1})_{2,7}=3/2$, so we are done.
            \item $E_8$: there are no components of multiplicity $1$ besides $C_0$, so there is nothing to check.
        \end{itemize}
    \end{itemize}
\end{proof}

\subsection{Section classes extend if and only if they are torsion}
Using the description of $\MW(\pi^\circ)$ in terms of monodromy coinvariants, we show that a section over $\partial D$ extends to a section over $D$ if and only if it is torsion in $\MW(\pi^\circ)$.
\begin{lemma}\label{lem: MW(pi^circ)tor extends}
    Let $\pi\colon X\to D$ be an elliptic fibration over a disk with one singular fiber at the origin.
   The image of the restriction map $\MW(\pi)\to \MW(\pi^\circ)$ equals $\MW(\pi^\circ)_\text{tor}$.
\end{lemma}
\begin{proof}
    By Proposition \ref{prop: MW(pi) if only one singularity}, $\MW(\pi)\cong G(X_0)$.
    By Corollary \ref{cor: MW pi is finite}, $\MW(\pi)$ is torsion. Thus, the image of the inclusion must be contained in the torsion part of $\MW(\pi^\circ)$. Furthermore, Lemma \ref{lem: component group is monodromy coinvariants} implies that $G(X_0)$ is isomorphic to the torsion part of the monodromy coinvariants, i.e.\ $G(X_0)\cong (\Z^2/(f_*-I)\Z^2)_\text{tor}$. Thus, the injection $\MW(\pi)\xhookrightarrow{} \MW(\pi^\circ)$ given by Proposition \ref{prop: isotopic on circle implies isotopic on disk} becomes an injection ${(\Z^2/(f_*-I)\Z^2)_\text{tor} \xhookrightarrow{} \MW(\pi^\circ)_\text{tor}}$.
    By Lemma \ref{lem: sections bij to coinvariants monodromy}, the codomain is equal to the domain and since both are finite, the map $\MW(\pi)\to\MW(\pi^\circ)_\text{tor}$ must be a bijection.
\end{proof}

\subsection{Proof of \Cref{thm: extend sections locally}}

\begin{proof}[Proof of \Cref{thm: extend sections locally}]
    Injectivity of the restriction map $\MW(\pi)\to \MW(\pi^\circ)$ follows from Proposition \ref{prop: isotopic on circle implies isotopic on disk}.
    We distinguish cases based on whether $X_0$ is additive or multiplicative. 
    
    \begin{enumerate}[(1)]
        \item Assume $X_0$ is additive. By definition, $\pi_1(X_0)=0$. Thus, $H_1(X_0)\cong 0$ and to prove exactness of the sequence, it suffices to show surjectivity of $\MW(\pi)\to \MW(\pi^\circ)$. 
        By Corollary \ref{cor: additive implies MW boundary finite}, $\MW(\pi^\circ)$ is a finite group. By Lemma \ref{lem: MW(pi^circ)tor extends}, this implies $\MW(\pi)\to \MW(\pi^\circ)$ is surjective, so the sequence is exact.

        \item Assume $X_0$ is multiplicative. By Kodaira's classification, there exists a basis $\{\alpha, \beta\}$ of $H_1(X_t)$ such that the monodromy $f\in \Mod(X_t)$ acts on $H_1(X_t)$ as
        $f_*=\begin{pmatrix}
        1 & n\\
        0 & 1
        \end{pmatrix}$. 

        By Lemma \ref{lem: sections bij to coinvariants monodromy}, $\eta \colon \MW(\pi^\circ)\to H_1(\partial X)$ is an injection.
        The zero-section $\sigma_0|_{\partial D}\colon \partial D\to \partial X$ induces a splitting of the short exact sequence \eqref{eq: SES MW pi^circ and H_1partial X}. Hence,
        \[
            H_1(\partial X)= \eta\bigg(\MW(\pi^\circ)\bigg)\oplus (\sigma_0|_{\partial D})_* \bigg(H_1(\partial D)\bigg).
        \]
        Also by Lemma \ref{lem: sections bij to coinvariants monodromy}, it follows 
        \[
            \MW(\pi^\circ)\cong \eta(\MW(\pi^\circ))\cong\frac{H_1(X_t)}{(I-f_*)H_1(X_t)} \cong (\Z/n\Z)\oplus\Z,
        \]
        with the $(\Z/n\Z)$-factor generated by $\alpha$ and the $\Z$-factor generated by $\beta$. (See Figure \ref{fig:explicit degeneration I_2}). Thus, 
        \[
            H_1(\partial X)=(\Z/n\Z)\langle \alpha\rangle \oplus \Z\langle\beta\rangle \oplus \Z\langle \gamma \rangle,
        \]
        where $\gamma=\sigma_0(\partial D)$.
        Restricting the codomain to the image of $\eta$ gives an isomorphism
        \[
            \eta\colon \MW(\pi^\circ)\to (\Z/n\Z)\langle \alpha\rangle \oplus \Z\langle\beta\rangle< H_1(\partial X).
        \]
        \begin{figure}
    \centering
    \includegraphics[width=0.4\linewidth]{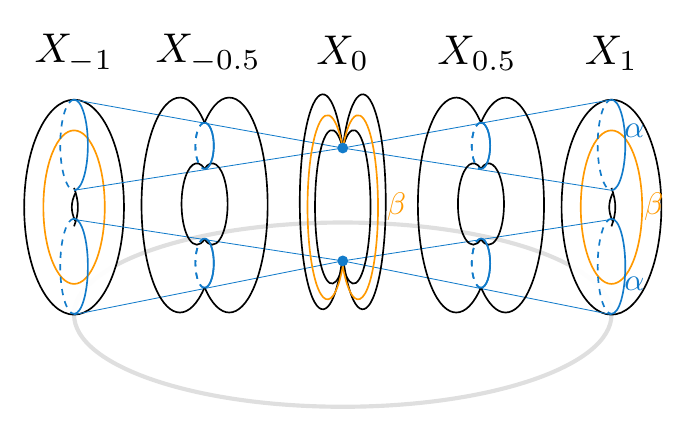}
    \caption{Explicit degeneration for a multiplicative fiber of type $I_2$.}
    \label{fig:explicit degeneration I_2}
\end{figure}
        Recall that we defined $\psi\colon \MW(\pi^\circ)\to H_1(X_0)$ as the composition of $\eta$ and the map $H_1(\partial X)\to H_1(X_0)$.
        By Picard-Lefschetz theory, $\alpha$ is the vanishing cycle of the fibration and thus its homology class vanishes in $H_1(X_0)$. The loop $\beta$ is mapped to a loop on $X_0$ going around the cycle of spheres precisely once, so $H_1(X_0)\cong \Z\langle \beta\rangle$.
        Therefore, the restriction of the map $H_1(\partial X)\to H_1(X_0)$ to $ (\Z/n\Z)\langle \alpha\rangle \oplus \Z\langle\beta\rangle$ takes $\alpha\mapsto0$, $\beta\mapsto \beta$. Therefore, the kernel of this restriction is $(\Z/n\Z)\langle \alpha \rangle$.
        Hence, $\ker(\psi)=\eta^{-1}((\Z/n\Z)\langle\alpha\rangle)$.

        By Lemma \ref{lem: MW(pi^circ)tor extends}, a section $\sigma\colon D^\circ\to X^\circ$ extends to a section over $D$ up to isotopy if and only if it is torsion in $\MW(\pi^\circ)$. Since $\eta$ is injective, $[\sigma]\in \MW(\pi^\circ)$ is torsion if and only if $\eta([\sigma])\in (\Z/n\Z)\langle \alpha\rangle \oplus \Z\langle\beta\rangle$ is torsion. The last statement holds if and only if $\eta([\sigma])\in (\Z/n\Z)\langle \alpha\rangle$. Thus $\MW(\pi)=\ker(\psi)$.

        Surjectivity of $\psi$ follows since the isotopy class of the section $\sigma\colon \partial D\to \partial X$ defined in Example \ref{ex: section does not extend} is mapped to the generator $[\beta]\in H_1(X_0)$. 
    \end{enumerate}
    
\end{proof}

\section{The Smooth Shioda--Tate Isomorphism} \label{section: compute full MW}
In this section, we compute the Mordell--Weil group of an arbitrary elliptic fibration $\pi$ over closed base.
We start by showing that $\MW(\pi)$ can be viewed as a subgroup of $\MW(\pi^\circ)$.
\begin{proposition} \label{prop: global injectivity Mw pi into MW pi circ}
    Let $\pi\colon M\to B$ be an elliptic fibration. Then the natural map $\MW(\pi)\to \MW(\pi^\circ)$ given by restriction is injective.
\end{proposition}
\begin{proof}
    If $[\sigma]\in \MW(\pi)$ and $\sigma|_{B^\circ}$ is isotopic to $\sigma_0|_{B^\circ}$, applying \Cref{prop: isotopic on circle implies isotopic on disk} in a small neighborhood of each critical value gives a local isotopy between the two sections, which can be adjusted to coincide with the original isotopy over $B^\circ$, yielding a global isotopy between $\sigma$ and $\sigma_0$.
\end{proof}

\subsection{A Global Obstruction to Extending Sections}
\label{section: globally: relation between MW(pi) and MW(pi^circ)}
In order to relate $\MW(\pi)$ and $\MW(\pi^\circ)$, we extend the local obstruction from \Cref{thm: extend sections locally} to a global one.
Consider an elliptic fibration $\pi\colon M\to B$ with critical value set $\Sigma\subset B$. For $t\in B$, denote $M_t\coloneqq \pi^{-1}(t)$.
Let $N$ be the closure of $M-M^\circ$ in $M$.
Define the map $\Psi\colon\MW(\pi^\circ)\to \bigoplus_{t\in \Sigma} H_1(M_t)$, as the composition of the following maps, which we explain below.
\[
    \MW(\pi^\circ)\xrightarrow[]{\iota} H_2(M^\circ,\partial M^\circ) \xrightarrow[]{\partial} H_1(\partial M^\circ)\xrightarrow[]{\cong} H_1(\partial N) \xrightarrow[]{\mathrm{incl}_*} H_1(N) \xrightarrow[]{\cong} H_1(\bigcup_{t\in\Sigma} M_t) \xrightarrow[]{\cong} \bigoplus_{t\in\Sigma} H_1(M_t),
\]

The homomorphism $\iota$ takes the isotopy class of a section $\sigma$ to $[\sigma]-[\sigma_0|_{B^\circ}]\in H_2(M^\circ,\partial M^\circ)$. The map $\partial\colon H_2(M^\circ,\partial M^\circ) \to H_1(\partial M^\circ)$ is the boundary map in the long exact sequence in homology of the pair $(M^\circ,\partial M^\circ)$. It acts on a homology class in $H_2(M^\circ,\partial M^\circ)$ by restricting it to the boundary $\partial M^\circ$. Observe that $\partial M^\circ$ and $\partial N$ are equal as sets, so the map $H_1(\partial M^\circ)\xrightarrow[]{\cong} H_1(\partial N)$ is induced by the identity. Inclusion $\partial N\xhookrightarrow{} N$ induces the map $H_1(\partial N)\xrightarrow[]{\mathrm{incl}_*} H_1(N)$. Since $N$ strong deformation-retracts to $\bigcup_{t\in \Sigma} M_t$, we obtain an isomorphism in homology $H_1(N) \xrightarrow[]{\cong} H_1(\bigcup_{t\in\Sigma} M_t)$. Finally, since the spaces $M_t$ are disjoint, homology of their union is isomorphic to the direct sum of their homology groups, yielding an isomorphism $H_1(\bigcup_{t\in\Sigma} M_t) \xrightarrow[]{\cong} \bigoplus_{t\in\Sigma} H_1(M_t)$.

Effectively, given an isotopy class $[\sigma]\in \MW(\pi^\circ)$, the map $\Psi$ restricts $\sigma$ to the boundary of a neighborhood of each singular fiber, obtaining a homology class in $\partial N$. Inclusion $\partial N\xhookrightarrow{} N$ gives a homology class in $N$, which deformation retracts to the union of the singular fibers.

For each singular value $t\in B$, let $D_t$ be a disk around $t$ containing no other critical values of $\pi$. Define the map $\psi_t\colon \MW(\pi|_{\pi^{-1}(D_t^\circ)})\to H_1(X_t)$ by $\psi_t([\sigma])=[\sigma|_{\partial D_t}]-[\sigma_0|_{\partial D_t}]$. This is the map $\psi$ from Section \ref{section: local picture}. We can use it to prove the following proposition.

\begin{proposition}[\textbf{A Global Obstruction to Extending Sections}]\label{prop: relate MW to MW(pi circ)}
    Let $\pi\colon M\to B$ be an elliptic fibration.
    The Mordell--Weil group $\MW(\pi)$ is the kernel of the map $\Psi\colon \MW(\pi^\circ)\to \bigoplus_{t\in\Sigma} H_1(M_t)$.
\end{proposition}
\begin{proof}
    A section $\sigma\colon B^\circ\to M^\circ$ extends to a section of $\pi$ if and only if it extends in a neighborhood of each fiber. We rewrite $\Psi$ in a way that lets us use Theorem \ref{thm: extend sections locally}.

    For $t\in \Sigma$, let $i_t\colon M_t\xhookrightarrow{} \pi^{-1}(\bar D_t)$ be the natural inclusion map.
    Homology of a union of disjoint union of spaces is naturally isomorphic to the direct sum of the homology groups of the spaces. Thus, there exists the following commutative diagram: 
    \[
    \begin{tikzcd}[row sep=large, column sep=huge]
    H_1(\partial N) \arrow[r, "\cong"] \arrow[d, "\incl_*"'] & 
    \bigoplus\limits_{t \in \Sigma} H_1(\pi^{-1}(\partial D_t)) \arrow[d, "\incl_*"] \\
    H_1(N) \arrow[r, "\cong"] \arrow[d, "\incl_*^{-1}"',"\cong"] & 
    \bigoplus\limits_{t \in \Sigma} H_1(\pi^{-1}(\overline{D}_t)) \arrow[d, "\cong"', "\bigoplus\limits_{t\in\Sigma}(i_t)_*^{-1}"] \\
    H_1\left(\bigcup\limits_{t \in \Sigma} M_t\right) \arrow[r, "\cong"] & 
    \bigoplus\limits_{t \in \Sigma} H_1(M_t),
    \end{tikzcd}
    \]
where the first vertical map is induced by inclusion $\partial N\xhookrightarrow{} N$ and the second vertical map is the inverse of the isomorphism induced by the inclusion $\bigcup_{t\in\Sigma} M_t\xhookrightarrow{} N$.
Hence, the composition of $\Psi$ with the projection $\bigoplus_{t \in \Sigma} H_1(M_t) \to H_1(M_t)$ is equal to the map $-\psi_t$. By \Cref{thm: extend sections locally}, $\sigma$ is isotopic to a section $\sigma_t\colon B^\circ\to M^\circ$ that is equal to $\sigma$ away from a neighborhood of $t$ and extends in a neighborhood of the critical value $t\in B$ if and only if $\psi_t([\sigma])=0$. 
Hence, $\sigma$ is isotopic over $B^\circ$ to a section $\sigma'\colon B^\circ \to M^\circ$, which extends to a global section $\sigma'\colon B\to M$ if and only if $\Psi([\sigma])=0$.

\end{proof}

\begin{remark}
    We will later see that the map $\MW(\pi^\circ)\to \bigoplus_{t\in \Sigma} H_1(M_t)$ is in fact surjective and the two Mordell--Weil groups fit in the following short exact sequence:
    \[
        0\to \MW(\pi)\to \MW(\pi^\circ)\to \bigoplus_{t\in \Sigma} H_1(M_t) \to 0.
    \]
\end{remark}

\begin{remark}
    In Proposition \ref{prop: relate MW to MW(pi circ)}, we can explicitly compute $\bigoplus_{t\in \Sigma} H_1(M_t)$. Indeed, the $M_t$ are simply connected if and only if they are additive fibers in Kodaira's classification. Each multiplicative fiber is a cycle of $n$ spheres and thus satisfies $H_1(M_t)\cong \Z$. Hence, if $m$ denotes the number of multiplicative fibers of $\pi$, it follows $\bigoplus_{t\in \Sigma} H_1(M_t)\cong \Z^m$. In particular, the Mordell--Weil groups of $\pi$ and $\pi^\circ$ differ by a finitely generated free group.
\end{remark}

\subsection{Proof of \Cref{thm: full MW}} \label{section: proof main thm}
We are finally ready to compute $\MW(\pi)$ in the general case where $\pi\colon M\to B$ is an elliptic fibration over a closed complex curve $B$.

Assume that the critical values of $\pi$ are $\Sigma=\{t_1,\dots,t_n\}\subset B$. Denote the irreducible components of the fiber $M_{t_i}$ by $C_{i,0},\dots, C_{i,r_i}$, where we assume that $C_{i,0}$ is the component intersecting the zero-section $\sigma_0$. These represent classes $[C_{i,j}]\in H_2(M)$. Let $c_{i,j}\in H^2(M)$ denote the Poincaré dual of $[C_{i,j}]$ in $M$. Let $e\in H^2(M)$ denote the Poincaré dual of the fiber class $F$ of $\pi$ and let $s\in H^2(M)$ denote the Poincaré dual of $[\sigma_0(B)]\in H_2(M)$.
Recalled that the \textit{trivial lattice} is defined by
\[
    \triv(M)\coloneqq\langle e,s, c_{i,j}: j>0 \rangle\subset H^2(M).
\]
We are now ready to prove \Cref{thm: full MW}.
\begin{proof}[Proof of \Cref{thm: full MW}]
    \step \label{step: rewrite LES}
    \textit{Construct an exact sequence relating $H^2(M), H^2(M^\circ)$ and $\bigoplus_{t\in \Sigma} H_1(M_t)$.}
    Consider the long exact sequence in cohomology corresponding to the pair $(M,M^\circ)$.
    \begin{equation}\label{eq: LES of pair}
        \dots \to H^1(M^\circ) \to H^2(M,M^\circ) \xrightarrow[]{j^*} H^2(M) \to H^2(M^\circ) \xrightarrow[]{\delta} H^3(M,M^\circ) \xrightarrow[]{j^*} H^3(M) \to \dots
    \end{equation}
    Let $N$ denote the closure of $M-M^\circ$. Observe that $N$ strong deformation retracts to $\bigcup_{t\in\Sigma} M_t$, so inclusion $\incl\colon \bigcup_{t\in\Sigma} M_t \xhookrightarrow{} N$ induces an isomorphism in homology.
    Let $l\colon (N,\partial N)\xhookrightarrow{} (M,M^\circ)$ be inclusion. By excision, $l$ induces an isomorphism $l^*\colon H^k(M,M^\circ)\to H^k(N,\partial N)$ for $k=0,1,2,3,4$.
    
    Define $\tau_k$ as the composition
    \[
    \begin{tikzcd}[column sep=large]
        H^k(M,M^\circ)\arrow[r,"l^*", "\cong"']\arrow[rrrr, bend right=10, "\tau_k","\cong"'] & H^k(N,\partial N) \arrow[r,"\PD_{[N, \partial N]}", "\cong"'] & H_{4-k}(N) \arrow[r, "(\incl_*)^{-1}", "\cong"'] & H_{4-k}\left(\bigcup\limits_{t\in \Sigma} M_t\right)\arrow[r, "\cong"'] &\bigoplus\limits_{t\in\Sigma} H_{4-k}(M_t),
    \end{tikzcd}
    \]
    where the last map is given by the fact that $M_t$ are disjoint sets.
    Using the isomorphisms $\tau_k$ and Poincaré duality in $M$, the long exact sequence (\ref{eq: LES of pair}) gives
    \begin{equation}\label{eq: LES after duality}
        \to H^1(M^\circ) \to \bigoplus_{t\in\Sigma} H_2(M_t) \xrightarrow[]{j^*\circ \tau_2^{-1}} H^2(M) \to H^2(M^\circ) \xrightarrow[]{\tau_3\circ\delta} \bigoplus_{t\in\Sigma} H_1(M_t) \xrightarrow[]{\PD_M\circ j^* \circ \tau_3^{-1}} H_1(M) \to \dots
    \end{equation}

    \step \label{step: map oplus H_k(X_t) to H_k(X) is incl}
    \textit{The map $\PD_M\circ j^*\circ (\tau_{4-k})^{-1}\colon \bigoplus_{t\in \Sigma} H_k(M_t)\to H_k(M)$ is induced by the inclusion $i_t\colon M_t\xhookrightarrow{} M$ on each summand, for $k=0,1,2,3,4$.}
    By definition of $\tau_{4-k}$, the claim is equivalent to showing that the rectangle in the diagram below commutes.
    
        \begin{center}
            \begin{tikzcd}[column sep=large]
        H^{4-k}(M, M^\circ)\arrow[dddddddd, bend right=50, "\tau_{4-k}"'] \arrow[rrr, "j^*"] \arrow[dd, "l^*"', "\cong"] &&&
        H^{4-k}(M) \arrow[dddddddd, "\cong"', "\PD_{M}"] \\
        &&&\\
        H^{4-k}(N, \partial N) \arrow[dd, "\cong", "\PD_{(N, \partial N)}"'] &&& \\
        &(1)&&\\
        H_k(N) \arrow[ddddrrr, "l_*", bend left=20] &&& \\
        &&&\\
        H_k(\bigcup\limits_{t\in \Sigma} M_t) \arrow[dd,"\cong"] \arrow[uu, "\incl_*", "\cong"'] \arrow[ddrrr, "\incl_*", bend left=10] &(2)&& \\
        &(3)&&\\
        \bigoplus\limits_{t\in\Sigma} H_k(M_t) \arrow[rrr,"\bigoplus\limits_{t\in\Sigma} (i_t)_*"'] &&&
        H_k(M) 
        \end{tikzcd}
        \end{center}
        
        As a composition of inclusions, triangle (2) commutes. Clearly, triangle (3) also commutes. Subdiagram (1) is commutative if and only if for any $\alpha\in H^{4-k}(M,M^\circ)$, the following identity holds in $H_k(M)$
        \begin{equation} \label{eq: frown step 2}
            l_*(l^*\alpha\frown [N,\partial N])=j^*\alpha\frown [M].
        \end{equation}
        By naturality of the cap product, the left-hand side is equal to $\alpha\frown l_*[N,\partial N]$. Since $l_*[N,\partial N]=j_*[M]\in H_4(M,M^\circ)$, the left hand side of \eqref{eq: frown step 2} equals $\alpha\frown j_*[M]$, so \eqref{eq: frown step 2} holds for all $\alpha\in H^{4-k}(M,M^\circ)$.

    \step \label{step: map from H_1 is 0}
    \textit{The map $\PD_M\circ j^* \circ (\tau_3)^{-1}\colon \bigoplus_{t\in\Sigma} H_1(M_t) \to H_1(M)$ in the long exact sequence (\ref{eq: LES after duality}) is zero.}
    By Step \ref{step: map oplus H_k(X_t) to H_k(X) is incl}, this map is $\bigoplus_{t\in\Sigma} (i_t)_*$.
    The map $\pi_*\colon H_1(M)\xrightarrow[]{} H_1(B)$ is an isomorphism by \Cref{prop: H_1 of elliptic fibration}. Thus, it suffices to show that $\pi_*\circ \bigoplus_{t\in\Sigma} (i_t)_*\colon \bigoplus_{t\in\Sigma} H_1(M_t) \to H_1(B)$ is zero. 
    But for each $t\in \Sigma$, the map $\pi\circ i_t$ is constant on $M_t$, so $\pi_*\circ (i_t)_*$ is zero on $H_1(M_t)$. Hence,
    the map $\PD_M\circ j^* \circ (\tau_3)^{-1}\colon \bigoplus_{t\in\Sigma} H_1(M_t) \to H_1(M)$ is trivial.

    \step \label{step: image oplus H_2}
    \textit{The image of the map $j^*\circ (\tau_2)^{-1}\colon \bigoplus_{t\in\Sigma}H_2(M_t)\to H^2(M)$ is $\langle c_{i,j}: i=\{1,\dots,n\}, j\ge 0\rangle$}.
    Observe that $\bigoplus_{t\in\Sigma} H_2(M_t)$ is the free abelian group with generators $[C_{i,j}]$. 
    By Step \ref{step: map oplus H_k(X_t) to H_k(X) is incl}, the following diagram commutes:
    \begin{center}
        \begin{tikzcd}
            H^2(M,M^\circ)\arrow[dd,bend right=60,"\tau_2"']\arrow[d,"\cong"]\arrow[r,"j^*"] &H^2(M)\arrow[d,"\cong"',"\PD_M"] \\
            H_2(\bigcup\limits_{t\in \Sigma}M_t)\arrow[d,"\cong"] \arrow[r,"\mathrm{incl}_*"]&H_2(M)\\
            \bigoplus\limits_{t\in\Sigma}H_2(M_t) &
        \end{tikzcd}
    \end{center}

    Hence, $j^*\circ(\tau_2)^{-1}([C_{i,j}])=\PD_M^{-1}[C_{i,j}]=c_{i,j}$.
    Thus, the image of the map $j^*\circ \tau_2^{-1}$ is precisely $\langle c_{i,j} \rangle \subset H^2(M)$. 
    
    \step \label{step: ker and annihilator}
    \textit{Deduce a short exact sequence relating the annihilators and replace $H^2(M^\circ)^F$ with $\MW(\pi^\circ)$.}
    By Step \ref{step: map from H_1 is 0} and Step \ref{step: image oplus H_2}, the long exact sequence (\ref{eq: LES after duality}) becomes
    \[
        0\to \langle c_{i,j} : i=\{1,\dots,n\}, j\ge 0\rangle \to H^2(M) \to H^2(M^\circ) \xrightarrow[]{\tau_3\circ\delta}  \bigoplus_{t\in\Sigma} H_1(M_t)\to 0.
    \]
    The components $C_{i,j}$ do not intersect the generic fiber, so all $c_{i,j}$ annihilate $F$. Restricting the exact sequence to the annihilators of $F$ induces the following exact sequence
    \[
        0\to \langle c_{i,j} : i=\{1,\dots,n\}, j\ge 0\rangle \to H^2(M)^F \to H^2(M^\circ)^F \xrightarrow[]{\tau_3\circ\delta} \bigoplus_{t\in\Sigma} H_1(M_t) \to 0.
    \]
    Surjectivity follows since for any $x\in H^2(M^\circ)$, it holds $x-\langle x,F\rangle s|_{M^\circ}\in H^2(M^\circ)^F$, which maps to the same element in $\bigoplus_{t\in\Sigma} H_1(M_t)$ as $x$.
    This induces the short exact sequence
    \begin{equation} \label{eq: SES H^2 X circ}
        0 \to H^2(M)^F/\langle c_{i,j} \rangle \to H^2(M^\circ)^F \xrightarrow[]{\tau_3\circ\delta} \bigoplus_{t\in\Sigma} H_1(M_t)\to 0.
    \end{equation}
    By \Cref{prop: MW bundle no singularities}, there exists an isomorphism $\Phi\colon \MW(\pi^\circ)\xrightarrow[]{\cong} H^2(M^\circ)^F$, which takes the isotopy class of a section $\sigma$ to $\PD[\sigma]-\PD[\sigma_0]\in H^2(M^\circ)^e$. Thus, the short exact sequence (\ref{eq: SES H^2 X circ}) induces the short exact sequence
    \begin{equation} \label{eq: SES MW(pi circ) as annihilator}
        0 \to H^2(M)^F/\langle c_{i,j} \rangle \to \MW(\pi^\circ) \xrightarrow[]{\tau_3\circ\delta\circ\Phi} \bigoplus_{t\in\Sigma} H_1(M_t)\to 0.
    \end{equation}

    By \Cref{prop: relate MW to MW(pi circ)}, $\MW(\pi)$ is the kernel of the map $\Psi\colon \MW(\pi^\circ)\to \bigoplus_{t\in \Sigma} H_1(M_t)$. Thus, to compute $\MW(\pi)$, it suffices to show that the two maps between $\MW(\pi^\circ)$ and $\bigoplus_{t\in \Sigma} H_1(M_t)$ are equal.

    \step \label{step: composition psi}
    \textit{The maps $\Psi\colon \MW(\pi^\circ)\to \bigoplus_{t\in \Sigma} H_1(M_t)$ and $\tau_3\circ\delta\circ \Phi \colon \MW(\pi^\circ)\to \bigoplus_{t\in \Sigma} H_1(M_t)$ are equal.}
    By the definitions of $\Phi$ and $\tau_3$, this is equivalent to showing commutativity of the following diagram:
    \begin{center}
        \begin{tikzcd}[row sep=large, column sep=huge]
    H^2(M^\circ) \arrow[r,blue, "\delta"] & 
    H^3(M, M^\circ) \arrow[d, "\cong"', "l^*"] \arrow[dddd,blue, "\tau_3", bend left=50] \\
    H_2(M^\circ, \partial M^\circ) \arrow[u, "\PD_{(M^\circ,\partial M^\circ)}^{-1}"', "\cong"] & 
    H^3(N, \partial N) \arrow[d, "\cong"', "\PD_{(N,\partial N)}"] \\
    \MW(\pi^\circ)\arrow[u] \arrow[uu,blue, "\Phi", bend left=60] 
    \arrow[ddr, blue, "\Psi"']  
    &H_1(N) \arrow[d, "\cong"'] \\
    & 
    H_1\left(\bigcup\limits_{t\in\Sigma} M_t\right) \arrow[d, "\cong"'] \\
    & 
    \bigoplus\limits_{t\in\Sigma} H_1(M_t)
\end{tikzcd}
    \end{center}

    We expand this diagram by writing out the maps in the definition of $\Psi\colon \MW(\pi^\circ)\to \bigoplus_{t\in \Sigma} H_1(M_t)$.
    We add some intermediary maps, which allow us to use naturality of the long exact sequence of a pair and of Poincaré--Lefschetz duality. 
        
    \noindent\hspace*{-2cm}
        \begin{tikzcd}
            & & H^2(M^\circ)\arrow[ddddll,"\frown {[ M^\circ,\partial M^\circ]}"',"\cong"]\arrow[dd,"\mathrm{restr}"]\arrow[rr,"="] & &H^2(M^\circ)\arrow[dd,"\mathrm{restr}"] \arrow[rr,blue,"\delta"]&& H^3(M,M^\circ)\arrow[dd,"l^*"] \arrow[dddddd,blue,bend left=40, "\tau_3"]\\
            &&& {\color{green!50!black}(1)} &&{\color{green!50!black}(2)}& \\
            && H^2(\partial M^\circ)\arrow[dd,"\frown {[\partial M^\circ]}","\cong"'] \arrow[rr,"\mathrm{id}^*","="'] &&H^2(\partial N) \arrow[dd,"\frown {[\partial N]}","\cong"']\arrow[rr,"\delta"] && H^3(N,\partial N)\arrow[dd,"\frown {[N,\partial N]}","\cong"']\\
            & {\color{green!50!black}(3)}& &{\color{green!50!black}(4)} & & {\color{green!50!black}(5)}& \\
            H_2(M^\circ,\partial M^\circ) \arrow[rr,"\partial"] &&H_1(\partial M^\circ) \arrow[rr,"\mathrm{id}_*","\cong"'] && H_1(\partial N) \arrow[rr,"\mathrm{incl}_*"] &&H_1(N)\arrow[d,"\cong"]\\
            &&&{\color{green!50!black}(6)}&&& H_1(\bigcup\limits_{t\in \Sigma} M_t)\arrow[d,"\cong"]\\ \MW(\pi^\circ)\arrow[uu] \arrow[uuuuuurr,blue,bend left=90,"\Phi"]\arrow[rrrrrr,blue,"\Psi"]&&&&&&\bigoplus\limits_{t\in\Sigma}H_1(M_t)
        \end{tikzcd}

    The first row is the identity, followed by the boundary map $\delta\colon H^2(M^\circ)\to H^3(M,M^\circ)$ from the long exact sequence of the pair $(M,M^\circ)$.
    The second row consists of the map induced by the identity $\mathrm{id}\colon \partial M^\circ\to \partial N$ and the boundary map $\delta\colon H^2(\partial N)\to H^3(N,\partial N)$ from the long exact sequence in cohomology of the pair $(N,\partial N)$. The third row is the boundary map $\partial \colon H_2(M^\circ,\partial M^\circ)\to H_1(\partial M^\circ)$, followed by the map on $H_1$ induced by the identity $\mathrm{id}\colon \partial M^\circ\to \partial N$, followed by the map $H_1(\partial N)\to H_1(N)$ induced by inclusion.

    We consider each square in the last diagram and show that it (anti-)commutes. Square $(1)$ commutes,  since both horizontal maps are the identity. The horizontal maps in square $(2)$ are the boundary maps from the long exact sequences of the pairs $(M,M^\circ)$, respectively $(N,\partial N)$ and the vertical maps are induced by the inclusion $(N,\partial N)\xhookrightarrow{}(M,M^\circ)$. Thus, square $(2)$ commutes by naturality of the long exact sequence of a pair. Square $(3)$ commutes and square $(5)$ anticommutes by naturality of Poincaré--Lefschetz duality with respect to the long exact sequences in homology and cohomology of a pair (see \cite[Theorem 9.2]{BredonTopology}). Square $(4)$ anticommutes, by naturality of the cap product and the fact that 
    $\partial M^\circ$ and $\partial N$ have opposite orientation, so $[\partial M^\circ]=-[\partial N]$.
    Square $(6)$ commutes, by definition of $\Psi$.

    Since $(1),(2),(4),(5)$ commute up to sign, with precisely 2 of them reversing sign, it follows that going from $H^2(M^\circ)$ to $H_1(N)$ by first going right and then down is the same as going down and then to the right. Since squares $(3)$ and $(6)$ commute, this shows that $\Psi=\tau\circ \delta\circ\Phi$.
    
    \step \label{step: annihilator is quotient by s} \textit{Show that $H^2(M)^F\cong H^2(M)/\Z s$}, where $s=\PD_M[\sigma_0(B)]$.
    Pairing with the fiber class induces the following short exact sequence:
    \begin{equation} \label{eq: SES annihilator}
        0\to H^2(M)^F \to H^2(M) \xrightarrow[]{\frown F}\Z\to 0,
    \end{equation}
    and since $\langle s,F\rangle = 1$, the homomorphism $\varphi\colon \Z\to H^2(M), \varphi(k)=k s$ induces a splitting of (\ref{eq: SES annihilator}). Hence,
    \[
        H^2(M)\cong H^2(M)^F \oplus \Z s
    \]
    and thus
    \[
        H^2(M)^F \cong H^2(M)/\Z s.
    \]
    
    \step \textit{Conclusion.}
    By applying in order \Cref{prop: relate MW to MW(pi circ)}, Step \ref{step: composition psi}, \eqref{eq: SES MW(pi circ) as annihilator} and Step \ref{step: annihilator is quotient by s}, we conclude
    \[
        \MW(\pi) = \ker(\Psi) =\ker(\tau_3\circ\delta\circ\Phi)\cong H^2(X)^F/\langle c_{i,j}\rangle \cong H^2(X)/\triv(\pi).
    \]
    The isomorphism can be given explicitly as follows. Consider the commutative diagram
    \[
    \begin{tikzcd}
        0\arrow[r] & H^2(M)/\triv(\pi) \arrow[r] & H^2(M^\circ)^F \arrow[r] &\bigoplus\limits_{t\in\Sigma} H_1(M_t) \arrow[r] & 0 \\
        0\arrow[r] & \MW(\pi) \arrow[r]\arrow[u, "\Phi_1"] & \MW(\pi^\circ) \arrow[r] \arrow[u,"\Phi"] &\bigoplus\limits_{t\in\Sigma} H_1(M_t) \arrow[r] \arrow[u,"\id"] & 0,
    \end{tikzcd}
    \]
    where 
    \[
        \Phi([\sigma])=\PD_{(M^\circ,\partial M^\circ)}([\sigma])-\PD_{(M^\circ,\partial M^\circ)}([\sigma_0|_{B^\circ}])
    \]
    for any $[\sigma]\in \MW(\pi^\circ)$, and 
    \[
        \Phi_1([\sigma])= [\PD_M(\sigma)-\PD_M(\sigma_0)] \in H^2(M)/\triv(\pi)
    \] 
    for any $[\sigma]\in \MW(\pi)$. By \Cref{prop: MW bundle no singularities}, $\Phi$ is an isomorphism, so the short five-Lemma implies that $\Phi_1$ is also an isomorphism.
\end{proof}

\begin{proof}[Proof of \Cref{thm: nielsen realization}]
    \Cref{thm: holomorphic Mordell--Weil} and \Cref{thm: full MW} give the following commutative diagram
    \[
    \begin{tikzcd}
        \MWH(\pi) \arrow[r, hook]\arrow[d,"\cong"] & \MW(\pi)\arrow[d,"\cong"] \\
        \NS(M)/\triv(\pi) \arrow[r, hook] & H^2(M)/\triv(\pi).
    \end{tikzcd}
    \]
    Since $H^2(M)/\NS(M)$ is torsion-free, the bottom arrow induces the identity on the torsion part, yielding the conclusion.
\end{proof}

\begin{proof}[Proof of \Cref{thm: relation MW pi and pi circ complete}]
    By \Cref{prop: global injectivity Mw pi into MW pi circ}, $\MW(\pi)$ injects into $\MW(\pi^\circ)$. By \Cref{prop: relate MW to MW(pi circ)}, $\MW(\pi)=\ker(\Psi)$.
    Surjectivity of $\Psi$ follows immediately from exactness of the short exact sequence \eqref{eq: SES MW(pi circ) as annihilator} and the fact that $\Psi=\tau_3\circ\delta\circ\Phi$.
\end{proof}

\subsection{Applications} \label{subsection: applications}
We give some explicit applications of \Cref{thm: full MW}.
\begin{proof}[Proof of \Cref{cor: rat ell surf characterization}]
    If $M$ is not a rational elliptic surface, \cite[Corollary 1.10]{Morariu_NarrowMW} implies that the rank of $\MW(\pi)$ is strictly larger than that of $\MWH(\pi)$. 
    If $M$ is a rational elliptic surface, then $\NS(M)=H^2(M)$, so \Cref{thm: holomorphic Mordell--Weil} and \Cref{thm: full MW} yield the conclusion.
\end{proof}

We conclude by elaborating on Example \ref{ex: intro Kummer} and Example \ref{ex: intro Inose} from the introduction, thus providing a complete computation of the Mordell--Weil group of these K3 surfaces as a direct application of \Cref{thm: full MW}.
\begin{ex}[\textbf{Mordell--Weil group of a Kummer surface}]
    Consider two non-isogenous elliptic curves $E_1,E_2$ with hyperelliptic involutions $\iota_1,\iota_2$. Let $M\coloneqq \mathrm{Km}(E_1\times E_2)$ be the minimal resolution of $(E_1\times E_2)/(\iota_1\times \iota_2)$. This can be viewed as an elliptic surface $\pi\colon M\to \p^1$ with $4$ singular fibers of type $I_0^*$, by considering projection $p_1\colon E_1\times E_2 \to E_1$ and the diagram
    \[
    \begin{tikzcd}[column sep=huge]
        E_1\times E_2 \arrow[r,dashed,"/(\iota_1\times \iota_2)"]\arrow[d,"p_1"] & M\arrow[d,"\pi"]\\
        E_1 \arrow[r,dashed, "/\iota_1"] &\p^1,
    \end{tikzcd}
    \]
    where the horizontal arrows are rational maps.
    For $1\le i\le 4$, let $a_i$ and $b_i$ denote the fixed points of $\iota_1$, respectively $\iota_2$. A section of $p_1$ is $s\colon E_1\to E_1\times E_2$, given by $s(x)=(x,b_1)$. Assume that the zero-section $\sigma_0\colon \p^1\to M$ is induced by $s$. We construct an explicit presentation for $H_2(M;\Z)$. For more details see \cite[Chapter VIII, Section 5]{BPV_CompactComplexSurfaces} and \cite[Section 2.2]{TaorminaWendland_SymmetryKummer}.
    
    Let $A_i\coloneqq (E_1\times \{b_i\})/(\iota_1\times 1)$ and $B_i\coloneqq (\{a_i\}\times E_2)/(1\times \iota_2)$. The minimal resolution of $(E_1\times E_2)/(\iota_1\times \iota_2)$ is obtained by replacing a neighborhood of the image of $(a_i,b_j)\in E_1\times E_2$ by the normal bundle of a sphere $C_{i,j}$ of self-intersection $-2$. For $i=1,2$ let $\alpha_i$ and $\beta_i$ be simple loops, such that $\alpha_i$ are setwise invariant under the involution $\iota_1$ and $\beta_i$ are setwise invariant under $\iota_2$. Assume that $\alpha_1\cdot \alpha_2=1$ and $\beta_1\cdot \beta_2=1$ and the homology classes of $\alpha_i$ and $\beta_i$ are generators of $H_1(E_1)$ and $H_1(E_2)$ respectively.
    The product $\alpha_i\times \beta_i$ is a torus in $E_1\times E_2$ containing $4$ fixed points of $\iota_1\times \iota_2$. Let $D_{i,j}'$ denote the image of $\alpha_i\times \beta_j-\{\text{fixed points}\}$ in $M$. By construction, $D_{i,j}'$ is a $4$ times punctured sphere.
    Let $D_{i,j}$ be the closed surface in $M$ obtained by capping the four punctures of $D_{i,j}'$.
    
    As divisors, the singular fibers of $\pi$ are $2B_i+\sum_{j=1}^4 C_{i,j}$ for $1\le i\le 4$ and the zero-section is $\sigma_0=A_1$. 
    Let $F$ denote the fiber class of $\pi$. Thus,
    \[
        \PD_M(\triv(\pi))= \Z\langle A_1,F\rangle \oplus \Z\langle B_i: 1\le i \le 4\rangle \oplus \Z\langle C_{i,j}: 1\le i\le 4 \text{ and } 2\le j\le 4 \rangle .
    \]
    All fibers are equal in homology, yielding the equation
    \begin{equation}\label{eq: fiber class F}
        2B_i+\sum_{j=1}^4 C_{i,j} = F \text{ for } 1\le i\le 4.
    \end{equation}
    Let $G$ be the fiber class of the fibration we would obtain from the projection $E_1\times E_2\to E_2$. Again, fibers are equal in homology, yielding the equation
    \begin{equation}\label{eq: fiber class G}
        2A_j+\sum_{i=1}^4 C_{i,j} = G \text{ for } 1\le j\le 4.
    \end{equation}
    Summing the equations in \eqref{eq: fiber class F} over $i$ and the ones in \eqref{eq: fiber class G} over $j$ yields the equation
    \[
        2\sum_{i=1}^4 A_i + 2\sum_{i=1}^4 B_i + 2\sum_{i,j=1}^4 C_{i,j}=4F+4G.
    \]
    Since $H_2(M;\Z)$ is torsion-free, this simplifies to
    \begin{equation}\label{eq: saturation}
        \sum_{i=1}^4 A_i + \sum_{i=1}^4 B_i + \sum_{i,j=1}^4 C_{i,j}-2F-2G=0.
    \end{equation}
    Let $L$ be the abelian group generated by 
    \[ 
        A_i,B_i,C_{i,j}, D_{k,l},F,G \text{ for } 1\le i,j\le 4, 1\le k,l\le 2
    \]
    with relations \eqref{eq: fiber class F}, \eqref{eq: fiber class G} and \eqref{eq: saturation}. By construction, $L$ is a sublattice of $H_2(M;\Z)$. It can be checked that $L\cong \Z^{22}$ and together with the intersection form, $L$ is unimodular. Hence $L=H_2(M;\Z)$.

    By \Cref{thm: full MW} and Poincaré duality,
    \[
        \MW(\pi)\cong H^2(M)/\triv(\pi) \cong H_2(M)/\PD_M(\triv(\pi)).
    \]
    The explicit presentations above imply that $\MW(\pi)$ is isomorphic to the group generated by 
    \[
        A_2,A_3, D_{k,l} \text{ for } 1\le k,l\le 2 
    \]
    with relations 
    \[
        2A_2=2A_3=0.
    \]
    Hence, 
    \[
        \MW(\pi)\cong \Z^4\times \Z/2\times \Z/2.
    \]
    This highlights the difference to the holomorphic case, where $\MWH(\pi)\cong \Z/2\times \Z/2$.
\end{ex}

\begin{figure}[h]
    \centering
    \includegraphics[width=0.6\linewidth]{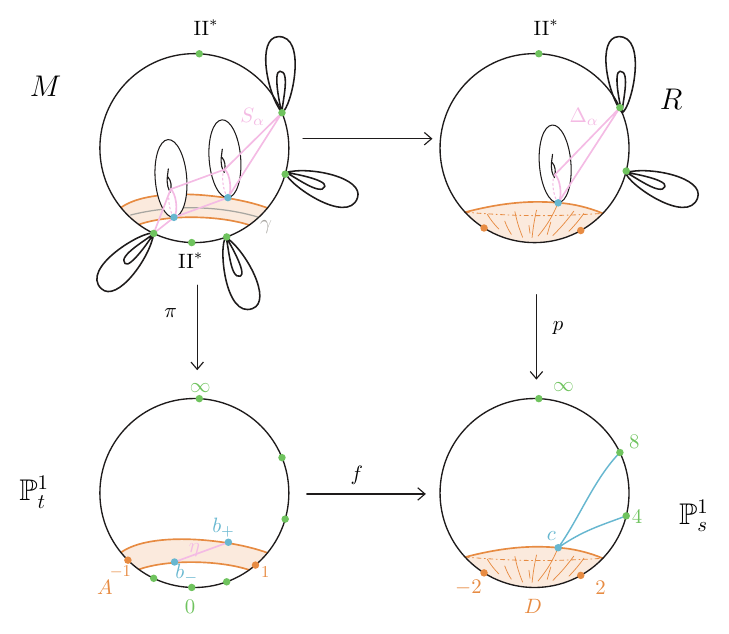}
    \caption{The Inose surface $\pi\colon M\to \p^1_t$ is obtained by pulling back the rational elliptic surface $p\colon R\to \p^1_s$ with $2$ nodal fibers and one fiber of type $II^*$. The branch values of $f\colon \p^1_t\to \p^1_s$ are regular points of $p$, so the singular fibers are doubled in $M$: there are four nodal fibers and two fibers of type $II^*$. The two pullbacks of a Lefschetz thimble of a nodal fiber of $p$ together with a cylinder in $\pi^{-1}(A)$ glue together to the embedded sphere $S_\alpha$.}
    \label{fig:inose construction}
\end{figure}

\begin{ex}[\textbf{Mordell--Weil group of an Inose surface}]\label{ex: inose surface}
    Let $p\colon R\to \p^1_s$ be the elliptic surface given as the minimal resolution of $y^2=x^3-3x+s-6$ with zero-section $\sigma_0$ at infinity. This is a rational elliptic surface with critical values $4,8,\infty$ and singular fibers of type $I_1,I_1$, respectively $II^*$.
    Consider the branched double cover $f\colon \p^1_t\to \p^1_s$ given by $f(t)=t+t^{-1}$. The branch values are $s=\pm 2$, neither of which is a singular value of $p$, so pulling the elliptic fibration $p$ back over $f$ yields an elliptically fibered K3 surface $\pi\colon M\to \p^1_t$ with two singular fibers of type $II^*$ and four singular fibers of type $I_1$, see \Cref{fig:inose construction}.

    Let $C_{0,1},\dots,C_0,8$ and $C_{\infty,1},\dots,C_{\infty,8}$ denote the irreducible components of the two fibers of type $II^*$ which do not intersect $\sigma_0$. Let $F$ denote a generic fiber. Then $C_{0,i}$, $C_{\infty,i}$, $F$ and $\sigma_0$ span the unimodular lattice
    \[
        E_8(-1)\oplus E_8(-1) \oplus 
        \begin{pmatrix}
            0 & 1 \\
            1 & -2
        \end{pmatrix}
    \]
    in $H_2(M)$.

    In order to determine the remaining $4$ generators, let $D$ be a disk in $P^1_s$ containing the branch values of $f$ but none of the critical values of $p$, as in \Cref{fig:inose construction}. The preimage $A\coloneqq f^{-1}(D)\subset \p^1_t$ is an annulus containing no critical values of $\pi$. Pick a point $c\in \partial D$ and fix paths in $\p^1_s-D$ from $c$ to the critical values $4,8,\infty$ of $p$ with mutually disjoint interiors, containing no other critical values of $p$. Let $\alpha,\beta\subset p^{-1}(c)$ be the vanishing cycles corresponding to the two nodal fibers of $p$. These have algebraic intersection number $1$. Let $\Delta_\alpha, \Delta_\beta$ be their vanishing thimbles in $R$ and assume that $\alpha,\beta,\Delta_\alpha,\Delta_\beta$ are disjoint from $\sigma_0$. Pick a trivialization of the torus-bundle $p\colon p^{-1}(D)\to D$ that preserves the zero-section. This induces a trivialization of $\pi$ over $A$. Let $\gamma$ be the central circle in the annulus $A$. Two further generators of $H_2(M)$ are given by the embedded tori $T_\alpha\coloneqq \alpha \times \gamma$ and $T_\beta\coloneqq \beta\times \gamma$.
    Let $f^{-1}(c)=\{b_-,b_+\}$. The pullback of $\Delta_\alpha$ to $M$ consists of two thimbles $\Delta_+$ and $\Delta_-$ with boundaries corresponding to the loop $\alpha$ in $M_{b_+}$, respectively $M_{b_-}$. The points $b_+$ and $b_-$ lie on different boundary components of $A$. Let $\eta$ be a properly embedded arc joining them, which intersects $\gamma$ transversely.
    The trivialization of $\pi$ over the annulus $A$ gives a cylinder $\alpha\times\eta$ in $M$. Together with this cylinder, the two thimbles glue to a continuously embedded sphere $S_\alpha$ in $M$, which can be smoothed. An analogous construction can be done for the preimage of $\Delta_\beta$, yielding a sphere $S_\beta$.

    By construction $T_\alpha,T_\beta,S_\alpha,S_\beta$ do not intersect the first $18$ generators we gave, so it suffices to show that these four homology classes generate a unimodular lattice in $H_2(M)$. Choosing a circle $\gamma'$ parallel to $\gamma$ in $A$ gives $[T_\alpha]= [\alpha\times \gamma']$, showing that the algebraic intersection numbers satisfy $T_\alpha\cdot T_\beta=0$, $T_\alpha\cdot T_\alpha=0$ and  $T_\beta\cdot T_\beta=0$.
    Observe that the algebraic intersection number of one of the tori with one of the spheres is equal to the algebraic intersection number of the vanishing cycles used in their respective construction. Hence, the intersection form on the lattice spanned by $T_\alpha,T_\beta,S_\alpha,S_\beta$ is
    \[
        \begin{pmatrix}
             & & & 1 \\
             & & -1& \\
             &1& *&*\\
            -1 & & *&*
        \end{pmatrix},
    \]
    so this lattice is unimodular.
    
    The homology classes $C_{0,1},\dots,C_{0,8}, C_{\infty,1},\dots, C_{\infty,8}, F, \sigma_0, T_\alpha,T_\beta,S_\alpha,S_\beta$ span a full rank unimodular sublattice of the unimodular lattice $H_2(M)$, so they span the entire $H_2(M)$. The trivial lattice is spanned by the Poincaré duals of $C_{0,1},\dots,C_{0,8}, C_{\infty,1},\dots, C_{\infty,8}, F, \sigma_0$, so \Cref{thm: full MW} implies
    \[
        \MW(\pi)\cong \Z^4
    \]
    and is generated by the Poincaré duals of $T_\alpha,T_\beta,S_\alpha,S_\beta$.
\end{ex}

\printbibliography
\end{document}